\documentclass[10pt]{amsart}

\usepackage[dvipsnames]{xcolor}
\usepackage{leftidx}
\usepackage{amsmath}
\usepackage{amscd}
\usepackage{amssymb}
\usepackage{rotating}
\usepackage{amsfonts}
\usepackage{amsthm}
\usepackage{verbatim}
\usepackage{stmaryrd}
\usepackage{mathrsfs}
\usepackage{color}
\usepackage[color]{xy}

\usepackage[colorlinks,linkcolor=blue, anchorcolor=blue, urlcolor=blue, citecolor=blue]{hyperref}
\usepackage{cleveref} 
 \usepackage[hyperpageref]{backref}

\usepackage[top=1in, bottom=1in, left=1.25in, right=1.25in]{geometry}

\input xy
\xyoption{all}

\allowdisplaybreaks[1]

\newcommand{\nn}{\nonumber}
\newcommand{\Spec}{\mathrm{Spec}}

\makeatletter
\newcommand{\holim@}[2]{%
  \vtop{\m@th\ialign{##\cr
    \hfil$#1\operator@font holim$\hfil\cr
    \noalign{\nointerlineskip\kern1.5\ex@}#2\cr
    \noalign{\nointerlineskip\kern-\ex@}\cr}}%
}
\newcommand{\holim}{%
  \mathop{\mathpalette\holim@{\leftarrowfill@\textstyle}}\nmlimits@
}
\newcommand{\hocolim@}[2]{%
  \vtop{\m@th\ialign{##\cr
    \hfil$#1\operator@font holim$\hfil\cr
    \noalign{\nointerlineskip\kern1.5\ex@}#2\cr
    \noalign{\nointerlineskip\kern-\ex@}\cr}}%
}
\newcommand{\hocolim}{%
  \mathop{\mathpalette\holim@{\rightarrowfill@\textstyle}}\nmlimits@
}
\makeatother

\newtheorem{theorem}{Theorem}[section]
\newtheorem{theorem/definition}{Theorem/Definition}[section]

\newtheorem{proposition}[theorem]{Proposition}
\newtheorem{lemma}[theorem]{Lemma}

\newtheorem{corollary}[theorem]{Corollary}

\theoremstyle{remark}
\newtheorem{remark}[theorem]{Remark}

\theoremstyle{definition}
 \newtheorem{example}[theorem]{Example}
\newtheorem{definition}[theorem]{Definition}

\newtheorem{problem}[theorem]{Problem}

\begin{document}
\begin{comment}
\title
{\large{\textbf{Rationality of dlog $\mathbb{A}^1$-zeta functions}}}
\author{\normalsize Xiaowen Hu}
\date{}
\maketitle
\end{comment}
\title{Rationality of $\mbox{dlog}$ $\mathbb{A}^1$-zeta functions}
\author{Xiaowen Hu}
\address{School of Sciences, Great Bay University, Dongguan, China}
\email{huxw06@gmail.com}
\dedicatory{In memory of Zhan-Wang Qian}

\begin{abstract}
For every smooth proper scheme over a finite field $\mathbb{F}_q$, Bilu, Ho, Srinivasan, Vogt, and Wickelgren introduced the dlog zeta function with coefficients in the Grothendieck-Witt ring $\mathrm{GW}(\mathbb{F}_q)$, enriching the dlog of the classical Weil zeta function with coefficients in $\mathbb{Z}$. They defined a notion of dlog rationality of such dlog zeta functions, which enriches the rationality of the Weil zeta function, and showed  the dlog rationality for simple cellular schemes. In this paper, we show that for any smooth proper schemes over $\mathbb{F}_q$, the dlog zeta function is rational, but not necessarily dlog rational. 
%Our proof is based on the rationality of the usual zeta functions, and is essentially elementary.
\end{abstract}

\maketitle

\section{Introduction}
For a smooth proper scheme $X$ over $\mathbb{F}_q$, the celebrated Weil conjecture says, among other assertions, that the zeta function
\begin{align*}
\zeta_X=\exp\big(\sum_m \frac{|X(\mathbb{F}_{q^{m}})|}{m}t^{m}\big)
\end{align*}
is a rational function of $t$. This is proved by Dwork \cite{Dwo60} and also by Grothendieck's school \cite{Gro64,Hou66,Del74,Del80}. In \cite{BHSVW22}, Bilu, Ho, Srinivasan, Vogt, and Wickelgren defined an analog of the Weil zeta function, by computing the trace of an endomorphism $\varphi$ of $X$ over a field $k$ with values in the endomorphism ring of the motivic sphere spectrum in the stable motivic homotopy category $\mathrm{SH}(k)$. A theorem of Morel \cite{Mor04} identifies this endomorphism ring with the Grothendieck–Witt ring $\mathrm{GW}(k)$ of the base field $k$.
By \cite[Chap.II. \S 2]{Lam05} for odd $q$, and (essentially) \cite[IV. Lemma 1.1]{MiH73} for even $q$ (see also \cite[Lemma 3.9]{Mor12}), we have
\begin{align*}
\mathrm{GW}(\mathbb{F}_q)=\begin{cases}
\mathbb{Z}\langle 1\rangle \cong \mathbb{Z},& \mbox{if}\ q\ \mbox{is even}\\
\frac{\mathbb{Z}[\langle u\rangle]}{-\langle 1\rangle+\langle u\rangle^2,-2\langle 1\rangle+2\langle u\rangle},& \mbox{if}\ q\ \mbox{is odd}
\end{cases}.
\end{align*}
Here $u$ denotes a non-square in $\mathbb{F}_q^{\times}$, and $\langle u\rangle$ denotes the quadratic form $(x,y)\mapsto uxy$. Note that $\langle 1\rangle=1$ in the ring $\mathrm{GW}(\mathbb{F}_q)$. For smooth proper schemes $X$ over $\mathbb{F}_q$ and  the relative Frobenius $\varphi$, an application \cite[Theorem 8.9]{BHSVW22} of Hoyois' trace formula \cite[Corollary 1.10]{Hoy14} yields 
\begin{align}\label{eq-traceFormula-quadratic}
\mathrm{dlog}\zeta^{\mathbb{A}^1}_{X}=& 
\sum_{m=1}^{\infty} \bigg(\sum_{\begin{subarray}{c}i\mid m\\ i\ even\end{subarray}}\big(\frac{1}{i}\sum_{d\mid i}\mu(d)|X(\mathbb{F}_{q^{i/d}})|\big)\big((i-1)\langle 1\rangle+\langle u\rangle\big)\nn\\
& +\sum_{\begin{subarray}{c}i\mid m\\ i\ odd\end{subarray}}\big(\frac{1}{i}\sum_{d\mid i}\mu(d)|X(\mathbb{F}_{q^{i/d}})|\big)i\langle 1\rangle\bigg)t^{m-1},
\end{align}
where $\mathrm{dlog}:=\frac{\mathrm{d}}{\mathrm{d}t}\circ\log$, and we have omitted $\varphi$ in the subscript of $\zeta^{\mathbb{A}^1}_{X}$.
We can rewrite (\ref{eq-traceFormula-quadratic}) as
\begin{align*}
\mathrm{dlog}\zeta^{\mathbb{A}^1}_{X}
=\langle 1\rangle \mathrm{dlog}\zeta_X+(-\langle 1\rangle+\langle u\rangle)\sum_{m\in 2 \mathbb{N}}\Big(\sum_{\begin{subarray}{c}i\mid m\\ i\ even\end{subarray}}\sum_{d\mid i}\frac{\mu(d)}{i}|X(\mathbb{F}_{q^{i/d}})|\Big) t^{m-1}.
\end{align*}
The term additional to the Weil zeta function is absent when $q$ is even. So we assume that $q$ is odd in the rest of this paper.

 A priori, $\mathrm{dlog}\zeta^{\mathbb{A}^1}_{X}\in \mathrm{GW}(\mathbb{F}_q)[[t]]$.  Since $\mathrm{GW}(\mathbb{F}_q)$ has torsions, we cannot define the exponential of $\mathrm{dlog}\zeta^{\mathbb{A}^1}_{X}$ (see further discussions on this issue in \cite[Remark 6.3]{BHSVW22}). Instead, the authors of loc. cit. introduced the notion of \emph{dlog rationality}.

%In this paper, we will denote $\mathrm{dlog}\zeta^{\mathbb{A}^1}_{X,\varphi}$ by $\mathrm{dlog}\zeta^{\mathbb{A}^1}_{X}$ for brevity.

\begin{definition}[{\cite[Definition 1.4]{BHSVW22}}]
A power series $h(t)\in \mathrm{GW}(\mathbb{F}_q)[[t]]$ is \emph{$\mbox{dlog}$ rational} if there exists finitely many polynomials $P_j\in \mathrm{GW}(\mathbb{F}_q)[t]$, and $c_j\in \mathrm{GW}(\mathbb{F}_q)$, $1\leq j\leq N$, such that $P_j\equiv 1\mod (t)$ in $\mathrm{GW}(\mathbb{F}_q)[t]$, and
\begin{gather*}
  h(t)=\sum_{j=1}^N c_j\frac{P_j'(t)}{P_j(t)}.
\end{gather*}
\end{definition}
Then using the \emph{cellular $\mathbb{A}^1$-homology} of Morel-Sawant \cite{MoS20}, they showed that $\mathrm{dlog}\zeta^{\mathbb{A}^1}_X$ is dlog rational for any smooth proper scheme $X/\mathbb{F}_q$ which admits \emph{a simple cellular structure} (see \cite[Theorem 2 and Definition 3.5]{BHSVW22}).
Let us see an example beyond the simple cellular ones. We set $\varepsilon=-\langle 1\rangle+\langle u\rangle$. Then
\begin{equation}\label{eq-relation-varepsilon=u-1}
  \varepsilon^2=0,\ 2\varepsilon=0.
\end{equation}

\begin{comment}
\begin{example}[{\cite[Example 8.12]{BHSVW22}}]
Let $X=\Spec(\mathbb{F}_{q^2})$.
\begin{align*}
&\mathrm{dlog}\zeta^{\mathbb{A}^1}_X=\sum_{m\ \mathrm{even}}(\langle 1\rangle +\langle u\rangle)t^{m-1}\\
=&\frac{2t}{1-t^2}+\frac{\varepsilon t}{1-t^2}= (1+\varepsilon)\mathrm{dlog}\frac{1}{1-(1+\varepsilon) t}+\mathrm{dlog} \frac{1}{1+t}.
\end{align*}
\end{example}
\end{comment}

\begin{example}\label{example-ellipticCurve-dlogRational}
Since
\begin{align*}
\mathrm{dlog}\frac{1}{1-t-t^2}\equiv \frac{1-t}{1-t^3}
\equiv  \frac{t+t^3}{1-t^6}+\frac{1+t^4}{1-t^6}\mod 2.
\end{align*}
and
\begin{align*}
\frac{t'+\big((t+t^2)t\big)'}{1-(t+t^2)^2}\equiv \frac{1+t^2}{(1+t+t^2)^2}\equiv \frac{1+t^4}{1-t^6}\mod 2,
\end{align*}
one has, by using (\ref{eq-dlog-1}) and (\ref{eq-dlog-2}) in Section \ref{sec:counterexample-dlog-rationality},
\begin{align*}
\mathrm{dlog}\frac{1}{1-t-t^2-\varepsilon t}- \mathrm{dlog}\frac{1}{1-t-t^2}=\frac{t+t^3}{1-t^6}\varepsilon .
\end{align*}
Let $E$ be the elliptic curve defined over $\mathbb{F}_5$ by $Y^2=X^3+X+1$. Then $a=1+5- \# E(\mathbb{F}_5)=-3$.
By Proposition \ref{prop-A1Zeta-ellipticCurve} (see also Example \ref{example-extensions-low-degrees}), this implies
\begin{align*}
\mathrm{dlog} \zeta^{\mathbb{A}^1}_E=\mathrm{dlog} \zeta_E+(-\langle 1\rangle+\langle u\rangle)\frac{t+t^3}{1-t^6},
\end{align*}
which is dlog rational as we have seen. We will show in Appendix \ref{sec:nonSimpleCellular-ellipticCurve} that elliptic curves over $\mathbb{F}_q$ are not cellular.
\end{example}

In this paper, we address the following problem.
\begin{problem}\label{question-rationalQuadZeta}
Let $X$ be a smooth proper scheme over $\mathbb{F}_q$. Are the following assertions true?
\begin{enumerate}
  \item[(i)] $\mathrm{dlog}\zeta^{\mathbb{A}^1}_{X,\varphi}$ is rational.
   \item[(ii)] $\mathrm{dlog}\zeta^{\mathbb{A}^1}_{X}$ is $\mathrm{dlog}$ rational.
 \end{enumerate} 
\end{problem}
Clearly (ii) implies (i). We will show:
\begin{theorem}\label{thm-main-rationality}
For any smooth proper scheme $X$ over $\mathbb{F}_q$, $\mathrm{dlog}\zeta^{\mathbb{A}^1}_{X,\varphi}$ is rational. Moreover, there exist examples such that $\mathrm{dlog}\zeta^{\mathbb{A}^1}_{X,\varphi}$ are not dlog rational.
\end{theorem}

By the Weil conjecture, there exist algebraic integers $\alpha_1,\dots,\alpha_r$ and $\beta_1,\dots,\beta_s$ such that
\begin{align*}
|X(\mathbb{F}_{q^n})|=\sum_{i=1}^r \alpha_i^n-\sum_{j=1}^s \beta_j^n.
\end{align*}
It follows that
\begin{align}\label{eq-dlogZeta-weilConj}
\mathrm{dlog}\zeta^{\mathbb{A}^1}_{X}-\mathrm{dlog}\zeta_{X}
=(-\langle 1\rangle+\langle u\rangle)\sum_{m\in2 \mathbb{N}}\Big(\sum_{j\mid m}(-1)^j\frac{\sum_{i=1}^r \alpha_i^j-\sum_{i=1}^s \beta_i^j}{j}\prod_{p\in \mathcal{P}(\frac{m}{j})}(1-\frac{1}{p})\Big)t^{m-1}.
\end{align}
Moreover, $\alpha_1,\dots,\alpha_r$ (resp. $\beta_1,\dots,\beta_s$) is the collection of roots (counting with multiplicity) of a monic polynomial in $\mathbb{Z}[X]$. This follows from the work of Dwork (see Koblitz' book \cite[Chapter V]{Kob84}), and also follows from the works of Deligne (use \cite[Théorème 1.6]{Del74} when $X$ is smooth projective, together with the purity result in \cite{Del80} when $X$ is smooth proper).

We will show that the above facts on $\zeta_X$ suffice to imply Theorem \ref{thm-main-rationality}.
 So we need to study the following problem in number theory.
Let $a_i\in \mathbb{Z}$ for $1\leq i\leq n$, and $\theta_1,\dots,\theta_n$ be the roots of the polynomial
\begin{align}\label{eq-polynomial-studyRootsPowerSum}
X^n-a_1X^{n-1}+a_2X^{n-2}+\dots+(-1)^na_n.
\end{align}
Denote $\mathbf{a}=(a_1,\dots,a_n)$. For natural numbers $m$, we define
\begin{align}\label{eq-F(m,a)}
F(m,\mathbf{a}):=\sum_{j\mid m}(-1)^j\frac{\sum_{i=1}^n\theta_i^j}{j}\prod_{p\in \mathcal{P}(\frac{m}{j})}(1-\frac{1}{p}).
\end{align}
\begin{comment}
Using Euler's function $\varphi$, we can rewrite
\begin{align*}
F(m,\mathbf{a}):=\frac{1}{m}\sum_{j\mid m}(-1)^j\varphi(\frac{m}{j})\sum_{i=1}^n\theta_i^j.
\end{align*}
\end{comment}
The rationality statement of Theorem \ref{thm-main-rationality} is then a consequence of:
\begin{theorem}\label{thm-periodicParity-intro}
Let $\mathbf{a}\in \mathbb{Z}^n$. Then $F(m,\mathbf{a}) \in \mathbb{Z}$ for $m\geq \mathbb{N}$, and $F(m,\mathbf{a}) \mod 2: \mathbb{N}\rightarrow \mathbb{Z}/2 \mathbb{Z}$, as a function of $m\in \mathbb{N}$, is eventually periodic. Moreover, if $a_n$ is odd, then $F(m,\mathbf{a}) \mod 2$ is periodic.
\end{theorem}
In fact, for Theorem \ref{thm-main-rationality} we need only to consider $\mathbf{a}$ with odd $a_n$, since $a_n$ is a power of $q$. But we insist on studying all $\mathbf{a}$ because this  requires only slightly more effort, and may have its own interest. We show also explicit bounds of the periods and the starting point of the periodicity, and the dependence on $\mathbf{a}\mod 4$. See Section \ref{sec:=rationality-and-periodicParity} for the precise statements. We explicitly compute the coefficient of $-\langle 1\rangle +\langle u\rangle $ in $\zeta^{\mathbb{A}^1}_E$ for elliptic curves $E$ (Proposition \ref{prop-A1Zeta-ellipticCurve}), and show the non-dlog rationality for certain elliptic curves by some elementary manipulations of rational functions.

%We do not address which are dlog rational...
The integrality of $F(m,\mathbf{a})$ is a direct consequence of the \emph{generalized Fermat's little theorem} (see Lemma \ref{lem-congruence-powerSum-1}), which is first shown by Schönemann \cite{Sch39}, and by Smyth \cite{Smy86} and Zarelua \cite{Zar06} in different ways. We give in Section \ref{sec:integrality-4periodicity} a (slightly) new proof, and some intermediate lemmas will be used again in this paper. The periodicity of $F(m,\mathbf{a}) \mod 2$ turns out to require kind of refined study of the generalized Fermat's little theorem. More precisely, putting
\begin{align*}
G(m,\mathbf{a}):=\frac{\sum_{i=1}^{n} \theta_i^{2m}-\sum_{i=1}^{n} \theta_i^{m}}{2^{v_2(m)+1}},
\end{align*}
then (see Proposition \ref{prop-F(m,a)=sumOfG(m,a)})
\begin{align*}
F(m,\mathbf{a})\equiv \sum_{k=1}^{v_2(m)}G(\frac{m}{2^k},\mathbf{a})\mod 2.
\end{align*}
The mod 2 periodicity of $G$ does not imply that of $F$ (see Example \ref{example-periodich-notImply-periodichhat}). We formulate this problem into a general form in Section \ref{sec:descendiblePeriodicFunc}, and find that we need exactly to show that  $G$ is a \emph{(essentially) descendible periodic function}, a notion that we coin in Section \ref{sec:descendiblePeriodicFunc}. We transform this notion into a convenient form. Finally, we complete the proof of  Theorem \ref{thm-periodicParity-intro} by proving similar properties of $G$ in local fields, similar to Zarelua's approach  to the generalized Fermat's little theorem. 

\section*{Notations}
In this paper, $\mu(n)$ denotes Möbius' function defined for $n\in\mathbb{N}$:
\begin{gather*}
\mu(n)=\begin{cases}
(-1)^r,& \mbox{if}\ n\ \mbox{has no repeated prime factors and $r$ is the number of prime factors of}\ n\\
0,& \mbox{if}\ n\ \mbox{has repeated prime factors}
\end{cases}.
\end{gather*}
For $n\in \mathbb{N}$, denote by $\mathcal{P}(n)$ the set of distinct prime factors of $n$. We will use the following formulae  without explicit references.
\begin{gather*}
\sum_{d\mid n}\mu(d)=\delta_{n,1},\\
\sum_{d\mid n}\frac{\mu(d)}{d}=\prod_{p\in \mathcal{P}(n)}(1-\frac{1}{p}).
\end{gather*}

\section*{Acknowledgment}
This paper arose out of Kirsten Wickelgren's lectures on GW-valued curve counting in the Chow-Witt summer school 2023 at BIMSA, which inspired me to learn about the GW-valued point counting problem in \cite{BHSVW22}. 
I am grateful to Nanjun Yang for organizing this  wonderful summer school, and to  Kirsten Wickelgren for helpful discussion on some related problems. I am grateful to Fabien Morel for explaining some details in \cite{MoS20} to me.  I am grateful to Mao Sheng, and Nanjun Yang again, for hosting my visit to YMSC and BIMSA this November, and interesting discussion that stimulates Appendix \ref{sec:nonSimpleCellular-ellipticCurve}. I  also thank  Peng Du, Yong Hu, Peng Sun, Yunhao Sun, Heng Xie, Jinxing Xu, and Lei Zhang for helpful discussions on various mathematics related to the topic of this paper.

%Besides the backgroud, the technique in this paper is essentially elementary, most of which I learned  in the high school. I owe a debt of gratitude to  my  coach Zhan-Wang Qian, and wish him rest in peace.

This work is supported by the start-up foundation of Great Bay University.

\section{Rationality and periodic parity}\label{sec:=rationality-and-periodicParity}
Recall the function $F(m,\mathbf{a})$ defined in (\ref{eq-F(m,a)}). We are going to show the following properties.
\begin{theorem}\label{thm-integrality}
For any $\mathbf{a}\in \mathbb{Z}^n$, $F(m,\mathbf{a})\in \mathbb{Z}$.
\end{theorem}

\begin{theorem}[Periodic parity]
\label{thm-parity-integralPolynomial}
For any $\mathbf{a}\in \mathbb{Z}^n$, there exists a natural number $N$, such that $2N\leq 2^{n+1}-2$,
\begin{enumerate}
  \item[(i)] $F(m,\mathbf{a})\mod 2$ is $2N$-periodic for $m\geq 2N+1$;
  \item[(ii)] if, in addition,  $a_n$ is odd, then $F(m,\mathbf{a})\mod 2$ is  $2N$-periodic for $m\in \mathbb{N}$.
\end{enumerate}
\end{theorem}

\begin{theorem}\label{thm-parityEquality-mod4}
$F(m,\mathbf{a})\mod 2$ is $4$-periodic in each component of $\mathbf{a}$. That means, for any $1\leq j\leq n$, if 
\begin{align*}
a'_i=\begin{cases}
a_i,& \mbox{if}\ i\neq j\\
a_i+4,& \mbox{if}\ i=j
\end{cases},
\end{align*}
then $F(m,\mathbf{a})\equiv F(m,\mathbf{a}')\mod 2$.
\end{theorem}
The explicit bounds in Theorem \ref{thm-parity-integralPolynomial} enable us to determine all values $F(m,\mathbf{a})\mod 2$ by computing $F(m,\mathbf{a})$ for $m\leq 2^{n+1}-2$ (resp. $m\leq 2^{n+2}-4$) in the case $a_n$ is odd (resp. $a_n$ is even) by brute force. Moreover, Theorem \ref{thm-parityEquality-mod4} enables us to display the results for all $\mathbf{a}$ of a given length.  We made a Macaulay2 program to do this. In the appendix, we give the results for $n\leq 3$.

\begin{example}\label{example-extensions-low-degrees}
By using Theorem \ref{thm-parityEquality-mod4},  we can compute $F(m,\mathbf{a})\mod 2$ by hand in some particular cases. For example, for $\mathbf{a}=(4,5)$, which is the case we will use in Example \ref{example-ellipticCurve-nonDlogRational}, we have $F(m,\mathbf{a})\equiv F(m,0,1)$.
But for $\mathbf{a}=(0,1)$, we have $\theta_i=\pm \sqrt{-1}$. Using Proposition \ref{prop-F(m,a)=sumOfG(m,a)}, one can then show 
\begin{align*}
F(m,\mathbf{a})\mod 2\equiv \begin{cases}
1,& \mbox{if}\ m\equiv 2\mod 4\\
0,& \mbox{otherwise}
\end{cases}.
\end{align*}
Another example is $\mathbf{a}=(-3,5)$ which is the case we used in Example \ref{example-ellipticCurve-dlogRational}. We have $F(m,\mathbf{a})\equiv F(m,1,1)$. For $\mathbf{a}=(1,1)$, we have $\theta_i=\zeta_6$ or $\zeta_6^5$, where $\zeta_6=\exp\frac{2\pi \sqrt{-1}}{6}$. Again by  Proposition \ref{prop-F(m,a)=sumOfG(m,a)} one gets
\begin{align*}
F(m,\mathbf{a})\mod 2\equiv \begin{cases}
1,& \mbox{if}\ m\equiv 2\ \mbox{or}\ 4 \mod 6\\
0,& \mbox{otherwise}
\end{cases}.
\end{align*}

\end{example}

\section{A counterexample to \texorpdfstring{$\mathrm{dlog}$}{dlog} rationality}\label{sec:counterexample-dlog-rationality}
Recall $\varepsilon=-\langle 1\rangle +\langle u\rangle$, where $u$ is a non-square element of $\mathbb{F}_q$. Then $\varepsilon$ is subject to the relations (\ref{eq-relation-varepsilon=u-1}).
If $f\in t \mathrm{GW}(\mathbb{F}_q)[[t]]\cong t(\mathbb{Z}\oplus \mathbb{Z}/2 \mathbb{Z})[[t]]$, we can write $f=g+\varepsilon h$, 
where $g\in t\mathbb{Z}[[t]]$ and $h\in t\frac{\mathbb{Z}}{2 \mathbb{Z}}[[t]]$. Then
\begin{align*}
f^n=(g+\varepsilon h)^n=g^n+n\varepsilon g^{n-1}h=\begin{cases}
g^n,& 2\mid n\\
g^n+\varepsilon g^{n-1}h,& 2\nmid n
\end{cases},
\end{align*}
\begin{align*}
\frac{1}{1-f}=\frac{1}{1-g}+\frac{\varepsilon h}{1-g^2},
\end{align*}
thus
\begin{align}\label{eq-dlog-1}
\mathrm{dlog} \frac{1}{1-f}=\frac{f'}{1-f}=\frac{g'+\varepsilon h'}{1-g}+\frac{\varepsilon g'h}{1-g^2}
=\frac{g'}{1-g}+\frac{h'+(gh)'}{1-g^2}\varepsilon ,
\end{align}
and
\begin{align}\label{eq-dlog-2}
\varepsilon \mathrm{dlog} \frac{1}{1-f}
=\frac{g'}{1-g}\varepsilon .
\end{align}
\begin{proposition}\label{prop-non-dlog-rationality}
Both $(-1+\langle u\rangle )\frac{t}{1-t^4}$ and $(-1+\langle u\rangle )\frac{t^3}{1-t^4}$
are not dlog rational.
\end{proposition}
\begin{proof}
Consider (\ref{eq-dlog-1}). The coefficients of odd powers of $t$ in $h'$, $(gh)'$, and $g^2$ are even. So the odd powers of $t$ in the coefficient of $\varepsilon $ in a $\mathrm{dlog}$ must arise from (\ref{eq-dlog-2}). Suppose that the leading terms of $g$ are $at+bt^2+ct^3+dt^4+\dots$. Then 
\begin{align*}
\mathrm{dlog} \frac{1}{1-g}=a+ \left(a^2+2 b\right)t+\left(a^3+3 a b+3 d\right)t^2+ \left(a^4+4 a^2 b+4 a d+2 b^2+4 e\right)t^3+\dots
\end{align*}
The coefficients of $t$ and $t^3$ have the same parity. So the assertion follows.
\end{proof}

\begin{proposition}\label{prop-A1Zeta-ellipticCurve}
Let $q$ be a power of an odd prime, and $E$ be an elliptic curve over $\mathbb{F}_q$. Let $a=1+q-\# E(\mathbb{F}_q)$. Then the coefficient of $\varepsilon t^{m-1}$ in $\mathrm{dlog}\zeta^{\mathbb{A}^1}_E$ is given in the following table.
\begin{table}[htbp!]
\centering
\begin{tabular}{|c|c|c|}
\hline
$(a,q)\mod 4$ & coefficient of $\varepsilon t^{m-1}$ in $\mathrm{dlog}\zeta^{\mathbb{A}^1}_E$\\
\hline
$(0,1)$ & =1 iff $m\in \{4k+2\}_{k\geq 0}$\\
$(0,3)$ &  0 \\
$(1,1)$ &  =1 iff $m\in\big\{6k+2,6k+4\}\big\}_{k\geq 0}$\\
$(1,3)$ &  =1 iff $m\in\{6k+4\}_{k\geq 0}$\\
$(2,1)$ & 0 \\
$(2,3)$ &  =1 iff $m\in\{4k+2\}_{k\geq0}$ \\
$(3,1)$ &  0\\
$(3,3)$ & =1 iff $m\in\{6k+2\}_{k\geq 0}$\\
\hline
\end{tabular}
\caption{Elliptic curves}
\label{tab:EllipticCurve}
\end{table}

\end{proposition}
\begin{proof}
By (\ref{eq-dlogZeta-weilConj}), the coefficient in consideration is equal to 
\begin{align*}
F(m,1)+F(m,q)-F(m,a,q).
\end{align*}
Then Table \ref{tab:EllipticCurve} follows from Tables \ref{tab:fMod2-n=1} and \ref{tab:fMod2-n=2}.
\end{proof}

\begin{example}\label{example-ellipticCurve-nonDlogRational}
Let $E$ be the elliptic curve defined over $\mathbb{F}_5$ by $Y^2=X^3+2X$. A direct counting yields $a=1+5- \# E(\mathbb{F}_5)=4$.
Then by Proposition \ref{prop-A1Zeta-ellipticCurve}, or using Example \ref{example-extensions-low-degrees}, we have
\begin{align*}
\mathrm{dlog} \zeta^{\mathbb{A}^1}_E=\mathrm{dlog} \zeta_E+(-1+\langle u\rangle)\frac{t}{1-t^4}.
\end{align*}
Hence by Proposition \ref{prop-non-dlog-rationality}, the statement (i) in Question \ref{question-rationalQuadZeta} in this case is not true. 
\end{example}
In the rest of this paper, I will show that the  statement (ii) in Question \ref{question-rationalQuadZeta} is true.

\section{Integrality and 4-periodic parity with respect to coefficients}\label{sec:integrality-4periodicity}
In this section, we show Theorem \ref{thm-integrality} and Theorem \ref{thm-parityEquality-mod4}, and make a first reduction of Theorem \ref{thm-parity-integralPolynomial} (see Proposition \ref{prop-F(m,a)=sumOfG(m,a)}).
Throughout this paper, we will use Newton's identity in the following form repeatedly:
\begin{align}\label{eq-Newton-identity}
\sum_{i=1}^n\theta_i^{m}=(-1)^m m \sum_{\begin{subarray}{c} r_1+2r_2+\dots+mr_m=m\\ r_1\geq 0,\dots,r_m\geq 0\end{subarray}}\frac{(r_1+r_2+\dots+r_m-1)!}{r_1!r_2!\cdots r_m!}\prod_{i=1}^m(-a_i)^{r_i}.
\end{align}

\begin{lemma}\label{lem-congruence-binomialCoefficient}
Let $p$ be a prime, and $a,b$ be natural numbers. Then
\begin{align*}
\binom{pa+pb}{pa}\equiv \binom{a+b}{a}\mod p.
\end{align*}
\end{lemma}
\begin{proof}
We compute
\begin{align*}
&\binom{pa+pb}{pa}- \binom{a+b}{a}=\frac{\prod_{k=1}^{pb}(pa+k)\cdot b!-\prod_{i=1}^b(a+i)\cdot (pb)!}{(pb)!b!}\\
=& \frac{\prod_{\begin{subarray}{c}1\leq k\leq pb\\ p\nmid k\end{subarray}}(pa+k)\cdot p^b\prod_{i=1}^b(a+i)\cdot b!
-\prod_{i=1}^b(a+i)\cdot \prod_{\begin{subarray}{c}1\leq k\leq pb\\ p\nmid k\end{subarray}}k \cdot p^b(b)!}{(pb)!b!}\\
=& \frac{p^b\prod_{i=1}^b(a+i)\Big(\prod_{\begin{subarray}{c}1\leq k\leq pb\\ p\nmid k\end{subarray}}(pa+k) 
-\prod_{\begin{subarray}{c}1\leq k\leq pb\\ p\nmid k\end{subarray}}k \Big)}{(pb)!}.
\end{align*}
Then the assertion follows from
\begin{align*}
v_p\big(p^b\prod_{i=1}^b(a+i)\big)\geq v_p(p^b b!)=v_p\big((pb)!\big),
\end{align*}
and 
\begin{align*}
v_p\Big(\prod_{\begin{subarray}{c}1\leq k\leq pb\\ p\nmid k\end{subarray}}(pa+k) 
-\prod_{\begin{subarray}{c}1\leq k\leq pb\\ p\nmid k\end{subarray}}k \Big)\geq 1.
\end{align*}
\end{proof}

\begin{lemma}\label{lem-integrality-factorial}
Let $p$ be a prime, and $m$ be a natural number. Let $r_i$ be natural numbers for $1\leq i\leq m$. Then
\begin{align*}
\frac{(pr_1+pr_2+\dots+pr_m-1)!}{(pr_1)!(pr_2)!\cdots (pr_m)!}-
p^{-1}\frac{(r_1+r_2+\dots+r_m-1)!}{r_1!r_2!\cdots r_m!}\in \mathbb{Z}_{(p)}.
\end{align*}

\end{lemma}
\begin{proof}
We do induction on $m$. The case $m=1$ is trivial. Suppose the assertion holds for $m-1$. Then
\begin{align*}
&\frac{(pr_1+pr_2+\dots+pr_m-1)!}{(pr_1)!(pr_2)!\cdots (pr_m)!}=\frac{\big(pr_1+pr_2+\dots+p_{m-2}+p(r_{m-1}+r_m)-1\big)!}{(pr_1)!(pr_2)!\cdots (pr_{m-2})!\big(p(r_{m-1}+r_m)\big)!}\cdot \binom{pr_{m-1}+pr_m}{pr_{m-1}}\\
\equiv & \left(p^{-1}\frac{(r_1+\dots+r_m-1)!}{r_1!\cdots r_{m-2}!(r_{m-1}+r_m)!} \right) \binom{pr_{m-1}+pr_m}{pr_{m-1}}\mod \mathbb{Z}_{(p)}.
\end{align*}
By Lemma \ref{lem-congruence-binomialCoefficient}, we have
\begin{align*}
  \binom{pr_{m-1}+pr_m}{pr_{m-1}}\equiv \binom{r_{m-1}+r_m}{r_{m-1}}\mod p.
\end{align*}
Hence the assertion holds for $m$.
\end{proof}

\begin{lemma}[{\cite{Sch39,Smy86,Zar06}}]
\label{lem-congruence-powerSum-1}
Let $p$ be a prime and $v\in \mathbb{N}$. Then
\begin{align*}
(-1)^{p^v}\sum_{i=1}^n\theta_i^{p^v}\equiv (-1)^{p^{v-1}}\sum_{i=1}^n\theta_i^{p^{v-1}}\mod p^v.
\end{align*}
\end{lemma}
\begin{proof}
Take $m=p^v$ in Newton's identity (\ref{eq-Newton-identity}). If there is some $r_k$ such that $p\nmid r_k$, then 
\begin{align*}
\frac{(r_1+r_2+\dots+r_m-1)!}{r_1!r_2!\cdots r_m!}\in \mathbb{Z},
\end{align*}
and thus
\begin{align*}
m \sum_{\begin{subarray}{c} r_1+2r_2+\dots+mr_m=m\\ r_1\geq 0,\dots,r_m\geq 0\end{subarray}}\frac{(r_1+r_2+\dots+r_m-1)!}{r_1!r_2!\cdots r_m!}\prod_{i=1}^m(-a_i)^{r_i}\equiv 0\mod p^v.
\end{align*}
Thus it suffices to show
\begin{gather*}
p^v\sum_{\begin{subarray}{c} r_1+2r_2+\dots+mr_m=p^{v-1}\\ r_1\geq 0,\dots,r_m\geq 0\end{subarray}}\frac{(pr_1+pr_2+\dots+pr_m-1)!}{(pr_1)!(pr_2)!\cdots (pr_m)!}\prod_{i=1}^m(-a_i)^{pr_i}\\
\equiv p^{v-1} \sum_{\begin{subarray}{c} r_1+2r_2+\dots+mr_m=p^{v-1}\\ r_1\geq 0,\dots,r_m\geq 0\end{subarray}}\frac{(r_1+r_2+\dots+r_m-1)!}{r_1!r_2!\cdots r_m!}\prod_{i=1}^m(-a_i)^{r_i}\mod p^v.
\end{gather*}
In fact, we are going to show, that for any $r_1,\dots,r_m\geq 0$ satisfying $r_1+2r_2+\dots+mr_m=p^{v-1}$, it holds that
\begin{gather}\label{eq-congruence-powerSum-0}
p^v\frac{(pr_1+pr_2+\dots+pr_m-1)!}{(pr_1)!(pr_2)!\cdots (pr_m)!}\prod_{i=1}^m(-a_i)^{pr_i}\\
\equiv p^{v-1}\frac{(r_1+r_2+\dots+r_m-1)!}{r_1!r_2!\cdots r_m!}\prod_{i=1}^m(-a_i)^{r_i}\mod p^v.\nn
\end{gather}
We have
\begin{align*}
v_p\big((p^ln)!\big)=&\sum_{k\geq 1}\lfloor \frac{p^ln}{p^k}\rfloor=(p^{l-1}+p^{l-1}+\dots+1)n+\sum_{k\geq l+1}\lfloor \frac{p^ln}{p^k}\rfloor=\frac{p^l-1}{p-1}n+v_p(n!),
\end{align*}
and
\begin{align*}
v_p\big((p^ln-1)!\big)=\frac{p^l-1}{p-1}n-l+v_p\big((n-1)!\big).
\end{align*}
Let $l=v_p\big(\mathrm{lcd}(r_1,\dots,r_m)\big)$, and $r_i=p^lr'_i$ for $1\leq i\leq m$. Then
\begin{align*}
&v_p\frac{(r_1+\dots+r_m-1)!}{r_1!\cdots r_m!}=v_p\frac{(p^l r'_1+\dots+p^l r'_m-1)!}{(p^lr'_1)!\cdots (p^lr'_m)!}\\
=& \frac{p^l-1}{p-1}\sum r'_i-l+v_p\big((\sum r'_i-1)!\big)
-\frac{p^l-1}{p-1}\sum r'_i+v_p\big(\sum r'_i)!\big)\\
=& -l+v_p\frac{(r'_1+\dots+r'_m-1)!}{r'_1!\dots r'_m !},
\end{align*}
and similarly
\begin{align}\label{eq-congruence-powerSum-1}
v_p\frac{(pr_1+\dots+pr_m-1)!}{(pr_1)!\cdots (pr_m)!}=-l-1+v_p\frac{(r'_1+\dots+r'_m-1)!}{r'_1!\dots r'_m !}\geq -l-1.
\end{align}
Since $p^l\geq l+1$, by Fermat's little theorem we have, for any integer $a$,
\begin{align*}
a^{p^{l+1}}-a^{p^l}=a^{p^l}(a^{p^{l+1}-p^l}-1)=a^{p^l}(a^{\varphi(p^{l+1})}-1)\equiv 0 \mod p^{l+1},
\end{align*}
and thus
\begin{align}\label{eq-congruence-powerSum-2}
\prod_{i=1}^m(-a_i)^{pr_i}\equiv \prod_{i=1}^m(-a_i)^{r_i}\mod p^{l+1}.
\end{align}
By (\ref{eq-congruence-powerSum-1}) and (\ref{eq-congruence-powerSum-2}) we have
\begin{gather*}
p^v\frac{(pr_1+pr_2+\dots+pr_m-1)!}{(pr_1)!(pr_2)!\cdots (pr_m)!}\prod_{i=1}^m(-a_i)^{pr_i}\\
\equiv 
p^v\frac{(pr_1+pr_2+\dots+pr_m-1)!}{(pr_1)!(pr_2)!\cdots (pr_m)!}\prod_{i=1}^m(-a_i)^{r_i}\mod p^v.
\end{gather*}
By Lemma \ref{lem-integrality-factorial}, we have
\begin{align*}
p^v\frac{(pr_1+pr_2+\dots+pr_m-1)!}{(pr_1)!(pr_2)!\cdots (pr_m)!}\equiv
p^{v-1}\frac{(r_1+r_2+\dots+r_m-1)!}{r_1!r_2!\cdots r_m!}\mod p^v.
\end{align*}
Hence (\ref{eq-congruence-powerSum-0}) follows.
\end{proof}

\begin{proof}[Proof of Theorem \ref{thm-integrality}]
Let $p$ be a prime, and $m=p^vm_0$ with $p\nmid m_0$. Then
\begin{align*}
F(m,\mathbf{a})=&\sum_{j\mid m_0}\sum_{k=0}^{v}\frac{(-1)^{p^k j}\sum_{i=1}^n\theta_i^j}{p^k j}\prod_{p'\in \mathcal{P}(\frac{m}{p^k j})}(1-\frac{1}{p'})\nn\\
=&\sum_{j\mid m_0}\sum_{k=0}^{v-1}\frac{(-1)^{p^k j}(p-1)\sum_{i=1}^n\theta_i^{p^k j}}{p^{k+1} j}\prod_{p'\in \mathcal{P}(\frac{m_0}{j})}(1-\frac{1}{p'})
+\frac{(-1)^{p^v j}\sum_{i=1}^n\theta_i^{p^vj}}{p^v j}\prod_{p'\in \mathcal{P}(\frac{m_0}{j})}(1-\frac{1}{p'})\nn\\
=&\sum_{j\mid m_0}\frac{\sum_{k=0}^{v-1}(-1)^{p^k j}(p-1)p^{v-k-1}\sum_{i=1}^n\theta_i^{p^k j}+(-1)^{p^v j}\sum_{i=1}^n\theta_i^{p^vj}}{p^{v} j}\prod_{p'\in \mathcal{P}(\frac{m_0}{j})}(1-\frac{1}{p'}).
\end{align*}
Note that
\begin{align}\label{eq-integrality-2}
&\frac{\sum_{k=0}^{v-1}(-1)^{p^k j}(p-1)p^{v-k-1}\sum_{i=1}^n\theta_i^{p^k j}+(-1)^{p^v j}\sum_{i=1}^n\theta_i^{p^vj}}{p^{v}}\nn\\
=&\sum_{k=1}^{v}\frac{(-1)^{p^k j}\sum_{i=1}^n\theta_i^{p^kj}-(-1)^{p^{k-1} j}\sum_{i=1}^n\theta_i^{p^{k-1}j}}{p^{k}}
+(-1)^{j}\sum_{i=1}^n\theta_i^{j}.
\end{align}
For each $j$ and $k \in \mathbb{N}$, applying Lemma \ref{lem-congruence-powerSum-1} to $\{\theta_i^{p^k j}\}_{i=1}^n$, which is also a complete set of conjugate algebraic integers, yields that 
\begin{align*}
\frac{(-1)^{p^k j}\sum_{i=1}^n\theta_i^{p^kj}-(-1)^{p^{k-1} j}\sum_{i=1}^n\theta_i^{p^{k-1}j}}{p^{k}}\in \mathbb{Z}.
\end{align*}
It follows that  $F(m,\mathbf{a})\in \mathbb{Z}_{(p)}$. This holds for each prime $p$, hence $F(m,\mathbf{a})\in \mathbb{Z}$.
\end{proof}

Define
\begin{align*}
G(m,\mathbf{a}):=\frac{\sum_{i=1}^{n} \theta_i^{2m}-\sum_{i=1}^{n} \theta_i^{m}}{2^{v_2(m)+1}}.
\end{align*}
By Lemma \ref{lem-congruence-powerSum-1}, $G(m,\mathbf{a})\in \mathbb{Z}$. 
A byproduct of the above proof is the following result.
\begin{proposition}\label{prop-F(m,a)=sumOfG(m,a)} 
Let  $m=2^vm_0$ with $2\nmid m_0$. Then
\begin{align*}
F(m,\mathbf{a})\equiv \sum_{k=1}^{v}G(\frac{m}{2^k},\mathbf{a})\mod 2.
\end{align*}
In particular, $F(m,\mathbf{a})$ is even if $m$ is odd. 
\end{proposition}
\begin{proof}
Consider the case $p=2$ in the above proof of Theorem \ref{thm-integrality}. Since $m_0$ is odd, we have
\begin{align*}
\prod_{p'\in \mathcal{P}(\frac{m_0}{j})}(1-\frac{1}{p'})\in 2\mathbb{Z}_{(2)}
\end{align*}
unless $m_0=1$. Hence $F(m,\mathbf{a})\equiv$ (\ref{eq-integrality-2}) with $j=1$ and $p=2$, which is equal to $\sum_{k=1}^{v}G(\frac{m}{2^k},\mathbf{a})$.
\end{proof}

As a consequence, Theorem \ref{thm-parityEquality-mod4} follows from the following result on $G(m,\mathbf{a})$.
\begin{theorem}\label{thm-parityEqualityOfG(m,e)-period4InE}
Let $m$ be an even natural number. Then $G(m,\mathbf{a})\mod 2$ is $4$-periodic in each component of $\mathbf{a}$.
\end{theorem}
The following lemma on the parity of $G(m,\mathbf{a})$ in the case $n=1$ is equivalent to Table \ref{tab:fMod2-n=1} in Example \ref{example-extensions-low-degrees}. Here we give a direct proof. 
\begin{lemma}\label{lem-tab:fMod2-n=1}
Let $m\in \mathbb{N}$. Then
\begin{align*}
G(m,a)\mod 2\equiv 
\begin{cases}
0,& \mbox{if}\ a\equiv 0\ \mbox{or}\ 1 \mod 4\\
\delta_{m,1}+\delta_{m,2},&\mbox{if}\ a\equiv 2\mod 4\\
\delta_{m,odd},&\mbox{if}\ a\equiv 3\mod 4
\end{cases}.
\end{align*}
\end{lemma}
\begin{proof}
Let $v=v_2(m)$, and $m=2^v m_0$. If $4\mid a$, then $4^m\mid a^m$. Since $2^{v+1}\geq v+2$, we have $4^{m}=2^{2^{v+1}m_0}\equiv 0\mod 2^{v+2}$.
So $2\mid G(m,a)$. 

If $a\equiv 2\mod 4$, let $a=2a_0$. We have $a^{m}=2^{2^vm_0}a_0^{m}$. Since $2^v\geq v+1$, and
\begin{align*}
2^{v}m_0< v+2\ \mbox{iff}\ (v,m_0)\in\{(0,1),(1,0)\},
\end{align*}
we obtain the result for $a\equiv 2\mod 4$.

Now let $a$ be odd, and write $a=\pm 1+4b$, with $b\in \mathbb{Z}$. By induction on $v$ we get 
\begin{align*}
 (\pm 1+4b)^{2^v}\equiv 1\mod 2^{v+2},\ \mbox{for}\ v\geq 1.
 \end{align*}
 So we get the result for $a\equiv 1\mod 4$. Moreover, $ (-1+4b)^{m_0}\equiv -1+4b\equiv -1\mod 4$.
 So 
 \begin{align*}
 G(m,a)\mod 2\equiv\begin{cases}
 0,& \mbox{if}\ 2\mid m\\
 1,& \mbox{if}\ 2\nmid m
 \end{cases}.
 \end{align*}
\end{proof}

\begin{proof}[Proof of Theorem \ref{thm-parityEqualityOfG(m,e)-period4InE}]
Let $m=2^vm_0$, where $2\nmid m_0$. By Newton's identity (\ref{eq-Newton-identity}), we have
\begin{align}
\sum_{i=1}^n\theta_i^{2m}-\sum_{i=1}^n\theta_i^{m}= &  2m \sum_{\begin{subarray}{c} r_1+2r_2+\dots+nr_n=2m\\ r_1\geq 0,\dots,r_n\geq 0\\ \exists\ odd\ r_j \end{subarray}}\frac{(r_1+r_2+\dots+r_n-1)!}{r_1!r_2!\cdots r_n!}\prod_{i=1}^n(-a_i)^{r_i}\label{eq-parityEqualityOfG(m,e)-period4InE-sum1}\\
&+2m \sum_{\begin{subarray}{c} r_1+2r_2+\dots+nr_n=m\\ r_1\geq 0,\dots,r_n\geq 0\end{subarray}}\frac{(2r_1+2r_2+\dots+2r_n-1)!}{(2r_1)!\cdots (2r_n)!}\prod_{i=1}^n(-a_i)^{2r_i}\label{eq-parityEqualityOfG(m,e)-period4InE-sum2}\\
&-m \sum_{\begin{subarray}{c} r_1+2r_2+\dots+nr_n=m\\ r_1\geq 0,\dots,r_n\geq 0\end{subarray}}\frac{(r_1+r_2+\dots+r_n-1)!}{r_1!r_2!\cdots r_n!}\prod_{i=1}^n(-a_i)^{r_i}.\label{eq-parityEqualityOfG(m,e)-period4InE-sum3}
\end{align}
 If there is some $r_k$ such that $2\nmid r_k$, then $\frac{(r_1+r_2+\dots+r_n-1)!}{r_1!r_2!\cdots r_n!}\in \mathbb{Z}_{(2)}$. Thus the sum (\ref{eq-parityEqualityOfG(m,e)-period4InE-sum1}) $\mod 2^{v+2}$ depends on the parities of $a_i$. We write the sum of (\ref{eq-parityEqualityOfG(m,e)-period4InE-sum2}) and (\ref{eq-parityEqualityOfG(m,e)-period4InE-sum3}) as
\begin{align*}
&2m \sum_{\begin{subarray}{c} r_1+2r_2+\dots+nr_n=m\\ r_1\geq 0,\dots,r_n\geq 0\end{subarray}}\frac{(2r_1+2r_2+\dots+2r_n-1)!}{(2r_1)!\cdots (2r_n)!}\left(\prod_{i=1}^n(-a_i)^{2r_i}-\prod_{i=1}^n(-a_i)^{r_i}\right)\\
&+2m\bigg(\sum_{\begin{subarray}{c} r_1+2r_2+\dots+nr_n=m\\ r_1\geq 0,\dots,r_n\geq 0\end{subarray}}\frac{(2r_1+2r_2+\dots+2r_n-1)!}{(2r_1)!\cdots (2r_n)!}\\
&-\frac{1}{2} \sum_{\begin{subarray}{c} r_1+2r_2+\dots+nr_n=m\\ r_1\geq 0,\dots,r_n\geq 0\end{subarray}}\frac{(r_1+r_2+\dots+r_n-1)!}{r_1!r_2!\cdots r_n!}\bigg)\prod_{i=1}^n(-a_i)^{r_i}.
\end{align*}
Let $l=v_p\big(\mathrm{lcd}(r_1,\dots,r_m)\big)$, and $r_i=2^l r'_i$. Then by (\ref{eq-congruence-powerSum-1}),
\begin{align*}
v_2\frac{(2r_1+\dots+2r_m-1)!}{(2r_1)!\cdots (2r_m)!}\geq -l-1.
\end{align*}
By Lemma \ref{lem-tab:fMod2-n=1},
\begin{align*}
\prod_{i=1}^m(-a_i)^{2r_i}-\prod_{i=1}^m(-a_i)^{r_i}
\equiv\begin{cases}
0\mod 2^{l+2},& \prod_{i=1}^m(-a_i)^{r'_i}\equiv 0\mod 4\\
0\mod 2^{l+2},& \prod_{i=1}^m(-a_i)^{r'_i}\equiv 1\mod 4\\
\delta_{l,0}2^{l+1}+\delta_{l,1}2^{l+1}\mod 2^{l+2},& \prod_{i=1}^m(-a_i)^{r'_i}\equiv 2\mod 4\\
\delta_{l,0}2^{l+1} \mod 2^{l+2},& \prod_{i=1}^m(-a_i)^{r'_i}\equiv 3\mod 4
\end{cases}.
\end{align*}
So
\begin{align*}
&\sum_{\begin{subarray}{c} r_1+2r_2+\dots+nr_n=m\\ r_1\geq 0,\dots,r_n\geq 0\end{subarray}}\frac{(2r_1+2r_2+\dots+2r_n-1)!}{(2r_1)!\cdots (2r_n)!}\Big(\prod_{i=1}^n(-a_i)^{2r_i}-\prod_{i=1}^n(-a_i)^{r_i}\Big)\in \mathbb{Z}_{(2)}
\end{align*}
and its parity in $\mathbb{Z}_{(2)}$ depends only on the mod 4 classes of the $a_i$'s.
By Lemma \ref{lem-integrality-factorial}, 
\begin{align*}
\sum_{\begin{subarray}{c} r_1+2r_2+\dots+nr_n=m\\ r_1\geq 0,\dots,r_n\geq 0\end{subarray}}\Big(\frac{(2r_1+2r_2+\dots+2r_n-1)!}{(2r_1)!\cdots (2r_n)!}-\frac{1}{2}\frac{(r_1+r_2+\dots+r_n-1)!}{r_1!r_2!\cdots r_n!}\Big)\prod_{i=1}^n(-a_i)^{r_i}\in \mathbb{Z}_{(2)}
\end{align*}
and its parity in $\mathbb{Z}_{(2)}$ depends only on the parities of the $a_i$'s. Hence the proof is completed.
\end{proof}

\section{Descendible periodic functions}\label{sec:descendiblePeriodicFunc}
Our final goal is to show the eventual periodicity of $F(m,\mathbf{a})\mod 2$. 
Proposition \ref{prop-F(m,a)=sumOfG(m,a)} motivates us to reduce the problem to a certain periodicity of $G(m,\mathbf{a})$.
\begin{definition}
Let $A$ be an abelian group, and $h:\mathbb{N}\rightarrow A$ a map. We define $\hat{h}:2\mathbb{N}\rightarrow A$ to be
\begin{align*}
\hat{h}(m)=\sum_{k=1}^{v_2(m)}h(\frac{m}{2^k}).
\end{align*}
\end{definition}
Then Proposition \ref{prop-F(m,a)=sumOfG(m,a)} says
\begin{align*}
F(m,\mathbf{a})\mod 2=\hat{G}(m,\mathbf{a})\mod 2,
\end{align*}
with both sides as functions of $m$, from $\mathbb{N}$ to $\mathbb{Z}/2 \mathbb{Z}$. The tables in the appendix give the values  of  $G(m,\mathbf{a})$ and $F(m,\mathbf{a})$ for $\mathbf{a}$ of length at most 3, in view of Theorem \ref{thm-parity-integralPolynomial} and \ref{thm-parityEquality-mod4} (see Example \ref{example-extensions-low-degrees}).

\begin{definition}
Let $A$ be a set. 
Let $h: \mathbb{N}\rightarrow A$ (resp. $2\mathbb{N}\rightarrow A$) be a map. We say that $h$ is \emph{periodic} if there exists $N\in \mathbb{N}$ such that 
\begin{align}\label{eq-periodicFunc}
h(x+N)=h(x)
\end{align}
for every $x\in  \mathbb{N}$, and in this case we say that $h$ is \emph{$N$-periodic}. Note that we do not demand that $N$ is a minimal period. If there exists $M\in \mathbb{N}$ (resp. $M\in 2 \mathbb{N}$) such that  (\ref{eq-periodicFunc}) holds for all $x\geq M$, then we say that $h$ is \emph{eventually ($N$-)periodic}.
\end{definition}

\begin{definition}
Let $A$ be a set. Let $h:A\rightarrow \mathbb{Z}/2 \mathbb{Z}$ be a $\mathbb{Z}/2 \mathbb{Z}$-valued function on $A$. We define the support of $h$ to be
\begin{align*}
\mathrm{Supp}(h):=\{x\in A: h(x)=1\}.
\end{align*}
\end{definition}

\begin{example}\label{example-periodich-notImply-periodichhat}
Let $h:\mathbb{N}\rightarrow \mathbb{Z}/2 \mathbb{Z}$ be a 4-periodic function. Then $\hat{h}$ is periodic if and only if $\mathrm{Supp}(h)$ is one of the following three cases:
\begin{align*}
\{4k+2\}_{k\geq 0},\ \{2k+1\}_{k\geq 0},\ \{4k+2\}_{k\geq 0}\cup\{2k+1\}_{k\geq 0}.
\end{align*}
We leave the proof to the reader.
\end{example}

Therefore periodicity of $h$ does not imply  the periodicity of $\hat{h}$.
For the latter, 
 %at least in the case that the minimal period of $\hat{h}$ is not larger than that of $h$, 
 $h$ has to satisfy a stronger condition that we investigate in the following.

\begin{definition}[Descendible periodic functions]\label{def-descendiblePeriodicFunc}
Let $N\in \mathbb{N}$. Let $h:\mathbb{N}\rightarrow \mathbb{Z}/2 \mathbb{Z}$  be a $2N$-periodic function. 
Let $S_0=\mathrm{Supp}(h)\cap [1,N]$, $S_1=\mathrm{Supp}(h)\cap [N+1,2N]$, and $S=S_0\sqcup S_1$. Let $\chi_{S_0}$ be the characteristic function of $S_0$, i.e.
\begin{align*}
\chi_{S_0}(x)=\begin{cases}
1,& \mbox{if}\ x\in S_0\\
0,& \mbox{if}\ x\in \mathbb{N}\backslash S_0
\end{cases}.
\end{align*}
We set
\begin{align*}
T:=& \mathrm{Supp}(\hat{\chi}_{S_0})\cap [1,2N]\\
=& \left\{x\in [1,2N]\cap 2 \mathbb{N}\Big| \#\big\{\frac{x}{2^k}| 1\leq k\leq v_2(x), h(\frac{x}{2^k})=1\big\}\ \mbox{is odd}\right\}.
\end{align*}
Then we say that $h$ is \emph{descendible ($2N$-)periodic} if the following conditions are satisfied, for all $x\in [1,N]\cap \mathbb{N}$:
\begin{enumerate}
  \item[(i)] if $2x\in T$ and $N+x\in T$, then $N+x\not\in S_1$;
  \item[(ii)] if $2x\in T$ and $N+x\not\in T$, then $N+x\in S_1$;
  \item[(iii)] if $2x\not\in T$ and $N+x\in T$, then $N+x\in S_1$;
  \item[(iv)] if $2x\not\in T$ and $N+x\not\in T$, then $N+x\not\in S_1$.
\end{enumerate}
\end{definition}
This definition is constructive. More precisely, given any subset $S_0$ of $[1,N]\cap \mathbb{N}$, the definition uniquely constructs a subset $S_1$ of $[N+1,2N]\cap \mathbb{N}$ such that $h$ is descendible $2N$-periodic. So we have:
\begin{lemma}
The descendible $2N$-periodic functions $\mathbb{N}\rightarrow \mathbb{Z}/2 \mathbb{Z}$ are 1-1 corresponding to subsets $S_0$ of $[1,N]\cap \mathbb{N}$.
\end{lemma}
For $x\in [1,N]\cap \mathbb{N}$, by the definition of $T$ we have
\begin{align*}
2x\in T\ \mbox{iff}\ 
(x\in S_0\ \&\ x\not\in T)\ \mbox{or}\ (x\not\in S_0\ \&\ x\in T).
\end{align*}
This provides a neater characterization of descendible periodic functions:
\begin{lemma}
Let $h:\mathbb{N}\rightarrow \mathbb{Z}/2 \mathbb{Z}$  be a $2N$-periodic function. Let $S$ and $T$ be the sets associated with $h$ as in Definition \ref{def-descendiblePeriodicFunc}. Then  $h$ is descendible $2N$-periodic if and only if
for all $x\in [1,N]\cap \mathbb{N}$,
\begin{align*}
\chi_{S}(x)+\chi_{S}(x+N)+\chi_{T}(x)+\chi_T(x+N)\equiv 0\mod 2.
\end{align*}
\end{lemma}

\begin{remark}\label{rmk-characterizationStrongPeriodicity}
Note that
\begin{gather*}
\chi_{S}(x)+\chi_{T}(x)\equiv \#\Big\{\frac{x}{2^k}\Big| 0\leq k\leq v_2(x),\ h(\frac{x}{2^k})=1 \Big\},\\
\chi_{S}(x+N)+\chi_{T}(x+N)\equiv \#\Big\{\frac{x+N}{2^k}\Big| 0\leq k\leq v_2(x+N),\ h(\frac{x+N}{2^k})=1 \Big\},
\end{gather*}
which give us a more direct characterization of descendible periodic functions.
\end{remark}
The following proposition is the motivation for defining descendible periodic functions.
\begin{proposition}\label{prop-periodic-and-descendiblePeriodic}
Let $h:\mathbb{N}\rightarrow \mathbb{Z}/2 \mathbb{Z}$  be a $2N$-periodic function. Then $\hat{h}:2 \mathbb{N}\rightarrow \mathbb{Z}/2 \mathbb{Z}$ is $2N$-periodic if and only if $h$ is descendible $2N$-periodic.
\end{proposition}
\begin{proof}
The “only if" direction  is clear. For example, $2x\in T\Longrightarrow \hat{h}(2x)=1 \Longrightarrow \hat{h}(2N+2x)=1$, if moreover $N+x\in T$, then $N+x\not\in S_1$. This is the case (i) in Definition \ref{def-descendiblePeriodicFunc}. The cases (ii)-(iv) are similar.

Now we show the “if" direction. For $x\in [1,N]\cap \mathbb{N}$, and $k\geq 0$, consider $\hat{h}\big((4k+2)N+2x\big)$. Since $\frac{(4k+2)N+2x}{2}=2kN+(N+x)$, we have
\begin{align*}
\hat{h}\big((4k+2)N+2x\big)=1\ \mbox{iff}\ (N+x\in S_1\ \&\ N+x\not\in T)\ \mbox{or}\ (N+x\not\in S_1\ \&\ N+x\in T)\ \mbox{iff}\ 2x\in T.
\end{align*}
Then for $k\geq 1$, consider $\hat{h}(4kN+2x)$. Since $\frac{4kN+2x}{2}=2kN+x$, we have
\begin{align*}
\hat{h}(4kN+2x)=1\ \mbox{iff}\ (x\in S_0\ \&\ x\not\in T)\ \mbox{or}\ (x\not\in S_0\ \&\ x\in T)\ \mbox{iff}\ 2x\in T.
\end{align*}
So the proof is completed.
\end{proof}

This proposition gives a complete characterization of the periodic functions $h:\mathbb{N}\rightarrow \mathbb{Z}/2 \mathbb{Z}$ such that $\hat{h}$ is periodic, since for any such $h$, we can find always a common even period for $h$ and $\hat{h}$. In the following, we study the eventually periodic functions.

\begin{definition}[Essentially descendible periodic functions]\label{def-essDescendiblePeriodicFunc}
Let $N\in \mathbb{N}$. Let $h:\mathbb{N}\rightarrow \mathbb{Z}/2 \mathbb{Z}$  be a eventually $2N$-periodic function. We say that $h$ is \emph{essentially descendible $2N$-periodic} if there exists a descendible $2N$-periodic function $g:\mathbb{N}\rightarrow \mathbb{Z}/2 \mathbb{Z}$ such that $h(i)=g(i)$ for $i\geq 2N+1$, and 
\begin{align}\label{eq-essDescendiblePeriodicFunc}
\sum_{k=0}^{2^km\leq 2N}h(2^km)\equiv \sum_{k=0}^{2^km\leq 2N}g(2^km)\mod 2.
\end{align}
\end{definition}

\begin{proposition}\label{prop-eventuallyPeriodic-and-essDescendiblePeriodic}
Let $h:\mathbb{N}\rightarrow \mathbb{Z}/2 \mathbb{Z}$  be an essentially descendible $2N$-periodic function. Then $\hat{h}:2 \mathbb{N}\rightarrow \mathbb{Z}/2 \mathbb{Z}$ is $2N$-periodic in the range $m\geq 2N+1$.
\end{proposition}
\begin{proof}
Let $g$ be a descendible $2N$-periodic function as in Definition \ref{def-essDescendiblePeriodicFunc}. Then from the definition of $\hat{h}$ and the assumptions on $g$ and $h$, we obtain $\hat{h}(i)=\hat{g}(i)$ for $i\geq 2N+1$. Hence the assertion follows from Proposition \ref{prop-periodic-and-descendiblePeriodic}.
\end{proof}

\begin{proposition}\label{prop-descendiblePeriodicit-sum}
If $h_i:\mathbb{N}\rightarrow \mathbb{Z}/2 \mathbb{Z}$  is  descendible (resp. essentially descendible) $2N_i$-periodic for $1\leq i\leq k$,  then $h=\sum_{i=1}^k h_i$ is  descendible (resp. essentially descendible) $\mathrm{lcm}(2N_1,\dots,2N_k)$-periodic.
\end{proposition}
\begin{proof}
By Proposition \ref{prop-periodic-and-descendiblePeriodic}, if $h_i:\mathbb{N}\rightarrow \mathbb{Z}/2 \mathbb{Z}$  is  descendible  $2N_i$-periodic, then $\hat{h}_i$ is $2N_i$-periodic, thus $\sum_{i=1}^k\hat{h}_i$ is $2\times\mathrm{lcm}(N_1,\dots,N_k)$-periodic, hence $\sum_{i=1}^k h_i$ is  descendible  $2\times\mathrm{lcm}(N_1,\dots,N_k)$-periodic.

If $h_i:\mathbb{N}\rightarrow \mathbb{Z}/2 \mathbb{Z}$  is essentially descendible  $2N_i$-periodic, let $g_i$ be a  descendible  $2N_i$-periodic as in Definition \ref{def-essDescendiblePeriodicFunc}. Then as above, $g=\sum_{i=1}^k g_i$ is $2\times\mathrm{lcm}(N_1,\dots,N_k)$-periodic, and we have $h(i)=g(i)$ for $i\geq 2\max(N_1,\dots,N_k)+1$. Then one can check (\ref{eq-essDescendiblePeriodicFunc}) holds for $N=\mathrm{lcm}(N_1,\dots,N_k)$.
\end{proof}

\section{Proof of the periodic parity}
In this section, we will show
\begin{theorem}\label{thm-periodicParity-G(m,a)}
For any $\mathbf{a}=(a_1,\dots,a_n)\in \mathbb{Z}^n$, $G(m,\mathbf{a})\mod 2$ is eventually descendible $2N$-periodic with $N\leq 2^{n}-1$. Moreover, if $a_n$ is odd, then $G(m,\mathbf{a})\mod 2$ is  descendible $2N$-periodic with $N\leq 2^{n}-1$.
\end{theorem}
Then Theorem \ref{thm-parity-integralPolynomial} follows as a consequence of Theorem \ref{thm-periodicParity-G(m,a)}, Proposition \ref{prop-F(m,a)=sumOfG(m,a)}, Proposition \ref{prop-periodic-and-descendiblePeriodic}, and Proposition \ref{prop-eventuallyPeriodic-and-essDescendiblePeriodic}. To show Theorem \ref{thm-periodicParity-G(m,a)}, we adapt the approach of Zarelua \cite{Zar06}.
%  to the generalized Fermat's little theorem. 
First we note that, by Proposition \ref{prop-descendiblePeriodicit-sum}, we can assume that (\ref{eq-polynomial-studyRootsPowerSum}) is irreducible over $\mathbb{Q}$. Then let $\mathbb{L}=\mathbb{Q}(\theta_1,\dots,\theta_n)$, and $\mathcal{O}_{\mathbb{L}}$ be the ring of integers in $\mathbb{L}$. Then
\begin{align*}
\sum_{i=1}^n\theta_i^m=\mathrm{Tr}_{\mathbb{Q}}^{\mathbb{L}}(\theta_1^m),\ \mbox{and}\
G(m,\mathbf{a})=\frac{\mathrm{Tr}_{\mathbb{Q}}^{\mathbb{L}}(\theta_1^{2m})-\mathrm{Tr}_{\mathbb{Q}}^{\mathbb{L}}(\theta_1^m)}{2^{v_2(m)+1}}.
\end{align*}
Since $\mathbb{L}/\mathbb{Q}$ is a Galois extension, there is a factorization
\begin{align*}
2\mathcal{O}_{\mathbb{L}}=\prod_{j=1}^r \mathfrak{P}_j^{e}
\end{align*}
of the ideal $2 \mathcal{O}_{\mathbb{L}}$ into prime ideals $\mathfrak{P}_j$ of $\mathcal{O}_{\mathbb{L}}$.  Let $\mathbb{L}_{\mathfrak{P}_j}$ be the local field associated with $\mathfrak{P}_j$. The degrees of extension of residue fields are the same for each $j$, and we denote it by $f$. So $n=ref$.
For $\xi\in \mathbb{L}$, we have (see e.g. \cite[II Cor.8.4]{Neu99})
\begin{align*}
\mathrm{Tr}_{\mathbb{Q}}^{\mathbb{L}}(\xi)=\sum_{j=1}^r \mathrm{Tr}_{\mathbb{Q}_2}^{\mathbb{L}_{\mathfrak{P}_j}}(\xi).
\end{align*}
Let 
\begin{align*}
G_{\mathfrak{P}_j}(m,\mathbf{a}):=\frac{\mathrm{Tr}_{\mathbb{Q}_2}^{\mathbb{L}_{\mathfrak{P}_j}}(\theta_1^{2m})-\mathrm{Tr}_{\mathbb{Q}_2}^{\mathbb{L}_{\mathfrak{P}_j}}(\theta_1^m)}{2^{v_2(m)+1}},\ \mbox{and}\
\hat{G}_{\mathfrak{P}_j}(m,\mathbf{a}):=\sum_{k=1}^{v_2(m)} G_{\mathfrak{P}_j}(\frac{m}{2^k},\mathbf{a}).
\end{align*}
Then
\begin{align}\label{eq-sum-G_i-localToGlobal}
G(m,\mathbf{a})=\sum_{j=1}^r G_{\mathfrak{P}_j}(m,\mathbf{a}),\ 
\hat{G}(m,\mathbf{a})=\sum_{i=1}^r\hat{G}_{\mathfrak{P}_j}(m,\mathbf{a}).
\end{align}

\begin{theorem}\label{thm-G(m)-essDescendiblePeriodic}
There exists a natural number  $N\leq 2^{ef}-1$ such that
$G_{\mathfrak{P}_j}(m,\mathbf{a})$ is essentially descendible $2N$-periodic. Moreover, if $a_n$ is odd, then $G_{\mathfrak{P}_j}(m,\mathbf{a})\mod 2$ is  descendible $2N$-periodic with $N\leq 2^{ef}-1$.
\end{theorem}
\begin{proof}[Proof of Theorem \ref{thm-periodicParity-G(m,a)}]
The assertion on descendible (resp. essentially descendible) periodicity follows from (\ref{eq-sum-G_i-localToGlobal}) and Proposition \ref{prop-descendiblePeriodicit-sum}. The upper bound on the periods follows by noting $2^{\sum_{j}n_j}-1\geq \prod_{j}(2^{n_j}-1)$.
\end{proof}

So we are reduced to show Theorem \ref{thm-G(m)-essDescendiblePeriodic}. The proof occupies the rest of this paper.
Let $L=\mathbb{L}_{\mathfrak{P}_j}$, and $\mathfrak{P}=\mathfrak{P}_j \mathcal{O}_L$. There is a unique sub-local field $K$ of $L$, such that  $L$ is totally ramified over $K$, and $K$ is unramified over $\mathbb{Q}_2$. Let $\mathfrak{p}=\mathfrak{P}\cap K$ be the maximal ideal of $\mathcal{O}_K$. Then $\mathcal{O}_L/{\mathfrak{P}}=\mathcal{O}_K/\mathfrak{p}$, $e=[L:K]$, and $f=[\mathcal{O}_K/\mathfrak{p}: \mathbb{F}_2]$. Let $\pi$ be a uniformizer of $L$. Then the minimal polynmial of $\pi$ over $K$
is an Eisenstein polynomial, and $\mathcal{O}_L$ is a free $\mathbb{Z}$-module with basis $\{1,\pi,\dots,\pi^{e-1}\}$. There is a $2^f-1$-th root of unity $w$ in $K$ such that $K=\mathbb{Q}_2(w)$.

For brevity of notations, in the following of this section we denote $\theta=\theta_1$, $\mathrm{Tr}=\mathrm{Tr}^L_{\mathbb{Q}_2}$, and $G(m)=G_{\mathfrak{P}_j}(m,\mathbf{a})$.
We define 
\begin{align*}
G^{L}_K(m):=\frac{\mathrm{Tr}_K^L(\theta^{2m})-\mathrm{Tr}_K^L(\theta^{m})}{2^{v_2(m)+1}}.
\end{align*}

%In this section, we denote $\mathrm{Tr}=\mathrm{Tr}_{K_{\mathfrak{p}}}^{L_{\mathfrak{P}}}$ for brevity of notations. 

\subsection{Odd \texorpdfstring{$a_n$}{an}}
That $a_n$ is odd is equivalent to  $\theta\not\in \mathfrak{P}$. Then $\theta=u+\xi$, where $\xi\in \mathfrak{P}$ and $u\in \mathcal{O}_K^{\times}$ satisfying $u^{2^f-1}=1$. Thus
\begin{align*}
&\mathrm{Tr}_K^L(\theta^{2m}-\theta^m)-e(u^{2m}-u^m)\\
=&\mathrm{Tr}_K^L\big((u+\xi)^{2m}-u^{2m}\big)-\mathrm{Tr}_K^L\big((u+\xi)^{m}-u^m\big).
\end{align*}
For $b\in K$ and $m\in \mathbb{N}$, we define
\begin{align*}
\tilde{G}_K^L(b,m):=\frac{\mathrm{Tr}_K^L\big((b+\xi)^{m}-b^m\big)}{2^{v_2(m)+1}}.
\end{align*}
 Let $e_i$ be the $i$-th elementary symmetric function of all conjugates of $\xi$ over $\mathbb{Q}$, where $1\leq i\leq e$. Since $\xi\in \mathfrak{P}$, we have $e_i\in \mathfrak{p}=2\mathcal{O}_K$. We let $e_i=2c_i$, with $c_i\in \mathcal{O}_K$.
\begin{lemma}\label{lem-crossingTerms-even-general}
Let $m\in \mathbb{N}$ and $b\in \mathcal{O}_K$. We expand
\begin{align*}
%\label{eq-minus-G(m)-mod2-oddTheta-general}
\frac{\mathrm{Tr}^L_K\big((b+\xi)^{m}-b^m\big)}{2^{v_2(m)+1}}
\end{align*}
as a polynomial in $c_i$'s, where $1\leq i\leq e$. Then the coefficient of every monomial of $c_i$'s with at least two distinct factors $c_i$ and $c_j$ (i.e. $i\neq j$) belongs to $2 \mathcal{O}_K$.
\end{lemma}
\begin{proof}
Recall Newton's identity (\ref{eq-Newton-identity}):
\begin{align}\label{eq-Newton-identity-1}
\sum_{i=1}^n\xi_i^{k}=(-1)^k k \sum_{\begin{subarray}{c} r_1+2r_2+\dots+er_{e}=k\\ r_1\geq 0,\dots,r_{e}\geq 0\end{subarray}}\frac{(r_1+r_2+\dots+r_{e}-1)!}{r_1!r_2!\cdots r_{e}!}\prod_{i=1}^{e}(-2c_i)^{r_i}.
\end{align}
Since
\begin{align*}
\frac{1}{2^{v_2(m)}}\binom{m}{k}k=\frac{m}{2^{v_2(m)}}\binom{m-1}{k-1}\in \mathbb{Z},
\end{align*}
it suffices to show
\begin{align}\label{eq-proof-crossingTerms-even-general}
\frac{(r_1+r_2+\dots+r_{2^{s-1}}-1)!}{r_1!r_2!\cdots r_{2^{s-1}}!}\prod_{i}2^{r_i}\in 4\mathbb{Z}_{(2)}
\end{align}
when there are at least two $r_i>0$. But 
\begin{align*}
\sum_{i}r_i-v_2(\sum_i r_i)\leq 1\ \mbox{iff}\ \sum_{i}r_i=1\ \mbox{or}\ 2.
\end{align*}
So we need only to consider the case $r_i=r_j=1$ for some $1\leq i<j\leq e$, and all other $r_k$'s are zero. In this case (\ref{eq-proof-crossingTerms-even-general}) is obvious. So the proof is completed.
\end{proof}

\begin{lemma}\label{lem-binomials-periodicity-general}
Let $s\in \mathbb{N}$. Then for $0\leq i\leq 2^s-1$ and $m\geq 0$, we have
\begin{align}\label{eq-binomials-periodicity-general}
\binom{m+2^s}{i}\equiv \binom{m}{i}\mod 2.
\end{align}
\end{lemma}
\begin{proof}By the formula
\begin{align*}
\binom{m}{i}=\binom{m-1}{i}+\binom{m-1}{i-1},
\end{align*}
we can do induction on $i$, and on $m$. Finally we arrive at the case $i=0$ and the case $m=0$. In the former case, both sides of (\ref{eq-binomials-periodicity-general}) are 0. In the latter case, we can assume $i>0$, then the left handside of (\ref{eq-binomials-periodicity-general}) is odd because $i\leq 2^s-1$.
\end{proof}

\begin{theorem}\label{thm-periodicity-G(m)-oddTheta-general}
Let $m\in \mathbb{N}$, and let $s\in \mathbb{N}$ such that $e\leq 2^{s-1}$. Then $\tilde{G}_K^L(b,m)\in \mathcal{O}_K$ for $b\in \mathcal{O}_K$, and
\begin{gather*}
      \tilde{G}_K^L(b,m+2^s)\equiv \tilde{G}_K^L(b,m)\mod 2 \mathcal{O}_K.
\end{gather*}
\end{theorem}
\begin{proof}
By Lemma \ref{lem-crossingTerms-even-general}, we need only to consider  the sum of powers of each single $c_i$ in $\tilde{G}_K^L(b,m)$, which by (\ref{eq-Newton-identity-1}) is equal to
\begin{align*}
\frac{1}{2^{v_2(m)+1}}
\sum_{r_i=1}^{\infty}\binom{m}{i r_i}b^{m-ir_i}(-1)^{i r_i}i(-2c_i)^{r_i}
=\frac{m}{2^{v_2(m)}}\sum_{r_i=1}^{\infty}\frac{2^{r_i}}{2r_i}\binom{m-1}{i r_i-1}b^{m-ir_i}(-1)^{i r_i}(-c_i)^{r_i}\in \mathcal{O}_K.
\end{align*}
This shows the first assertion. Note that $v_2\big(\frac{2^{r_i}}{2r_i}\big)<1$ exactly when $r_i=1$ or $2$. So we have
\begin{align}\label{eq-periodicity-G(m)-oddTheta-general}
&\frac{1}{2^{v_2(m)+1}}
\sum_{r_i=1}^{\infty}\binom{m}{i r_i}b^{m-ir_i}(-1)^{i r_i}i(-2c_i)^{r_i}\nn\\
\equiv& \frac{m}{2^{v_2(m)}}\left(\binom{m-1}{i-1}b^{m-i}(-1)^{i}(-c_i)
+\binom{m-1}{2i-1}b^{m-2i}c_i^{2}\right)\mod 2 \mathcal{O}_K.
%\equiv& u^m\left(\binom{m-1}{i-1}\frac{c_i}{u}+\binom{m-1}{2i-1}(\frac{c_i}{u})^2\right)\mod 2.
\end{align}
Then the second assertion follows from Lemma \ref{lem-binomials-periodicity-general}.
\end{proof}

\begin{comment}
 The following Lemma \ref{lem-congruence-totRam-general-1} is a special case of \cite[Lemma 5]{Zar06}. But the proof of the latter uses \cite[Lemma 4]{Zar06}, the proof of which I don't understand.
\begin{lemma}\label{lem-congruence-totRam-general-1}
Assume $a \in K$, $b\in \mathfrak{P}$, and $m\in \mathbb{N}$. Then
\begin{align*}
\mathrm{Tr}_K^L\big((a+b)^{m}-a^m\big)\in \mathfrak{p}^{v_2(m)+1}.
\end{align*}
\end{lemma}
\begin{proof}
Follows from the proof of Lemma \ref{lem-crossingTerms-even-general} and Theorem \ref{thm-periodicity-G(m)-oddTheta-general}.
\end{proof}

\begin{lemma}\label{lem-congruence-totRam-general-2}
Assume $a \in \mathcal{O}_K^{\times}$, $b\in \mathfrak{P}$, and $e\leq 2^{k-1}$. Then
\begin{align*}
\mathrm{Tr}_K^L\big((a+b)^{2^k}-a^{2^k}\big)\in \mathfrak{p}^{k+2}.
\end{align*}
\end{lemma}
\begin{proof}
This can be regarded as the case $m=0$ of the  following Theorem \ref{thm-periodicity-G(m)-oddTheta-general}.
\end{proof}
\end{comment}

By the first assertion of Theorem \ref{thm-periodicity-G(m)-oddTheta-general}, we have
\begin{align*}
\frac{\mathrm{Tr}_K^L(\theta^{2m}-\theta^m)}{2^{v_2(m)+1}}-\frac{e(u^{2m}-u^m)}{2^{v_2(m)+1}}
=2\tilde{G}_K^L(u,2m)-\tilde{G}_K^L(u,m)
\equiv  -\tilde{G}_K^L(u,m)\mod 2 \mathcal{O}_K.
\end{align*}
Then using $u^{2^f-1}=1$, the second assertion of Theorem \ref{thm-periodicity-G(m)-oddTheta-general}, and $G(m)=\mathrm{Tr}^{K}_{\mathbb{Q}_2}\circ G^L_K(m)$, we get
\begin{align*}
G\big(m+2^{s}(2^f-1)\big)\equiv G(m)\mod 2.
\end{align*}
We can take $s$ to be the smallest natural number such that $e\leq 2^{s-1}$. Then $G(m)$ is $2^{s}(2^f-1)$-periodic with $2^{s-1}(2^f-1)\leq 2^{ef}-1$. 
Then by Remark \ref{rmk-characterizationStrongPeriodicity}, the descendible periodicity of $G$ amounts to the following theorem.
\begin{theorem}\label{thm-descendiblePeriodicity-G(m)-oddTheta-general}
Suppose $e\leq 2^{s-1}$. Then for $1\leq m\leq 2^{s-1}(2^f-1)$,
\begin{align*}
\sum_{k=0}^{v_2(m)}\frac{\mathrm{Tr}\big((u+\xi)^{\frac{m}{2^k}}-u^{\frac{m}{2^k}}\big)}{2^{v_2(m)-k+1}}
\equiv\sum_{k=0}^{v_2\big(m+2^{s-1}(2^f-1)\big)}\frac{\mathrm{Tr}\big((u+\xi)^{\frac{m+2^{s-1}(2^f-1)}{2^k}}-u^{\frac{m+2^{s-1}(2^f-1)}{2^k}}\big)}{2^{v_2\big(m+2^{s-1}(2^f-1)\big)-k+1}}\mod 2.
\end{align*}
\end{theorem}

\begin{comment}
In particular, the case $m=2^{s-1}$ is equivalent to Lemma \ref{lem-congruence-totRam-general-2}.

Theorem \ref{thm-descendiblePeriodicity-G(m)-oddTheta-general} follows from Lemma \ref{lem-crossingTerms-even-general} and the following Lemma \ref{lem-binomials-sum-mod2Equality-general}. In the rest of this section, $\mathrm{Tr}:= \mathrm{Tr}^L_{\mathbb{Q}_2}$.
We take $u$ satisfying
\begin{align}\label{eq-fermatLittle-unit}
u^{m+2^s(2^f-1)}=u^m,\ e\leq 2^{s-1},
\end{align}
\end{comment}
By the proof of Theorem \ref{thm-periodicity-G(m)-oddTheta-general}, especially (\ref{eq-periodicity-G(m)-oddTheta-general}), we see that 
we need only to consider the sum of powers of each single $c_i$ in $\tilde{G}_K^L(u,l)$ (where $l$ runs over $\frac{m}{2^k}$ and $\frac{m+2^{s-1}(2^f-1)}{2^k}$). This sum is equal to
\begin{align*}
& \frac{l}{2^{v_2(l)}}\left(\binom{l-1}{i-1}u^{l-i}(-1)^{i}(-c_i)
+\binom{l-1}{2i-1}u^{l-2i}c_i^{2}\right)\\
\equiv& u^l\left(\binom{l-1}{i-1}\frac{c_i}{u}+\binom{l-1}{2i-1}(\frac{c_i}{u})^2\right)\mod 2 \mathcal{O}_K.
\end{align*}
Then Theorem \ref{thm-descendiblePeriodicity-G(m)-oddTheta-general} follows by applying the following Lemma \ref{lem-binomials-sum-mod2Equality-general} to $d_i:=\frac{c_i}{u}$.

\begin{lemma}\label{lem-a+a^2-trace-even}
For $\beta\in \mathcal{O}_K$, we have $\mathrm{Tr}^K_{\mathbb{Q}_2}(\beta+\beta^2)\in 2 \mathbb{Z}_2$.
\end{lemma}
\begin{proof}
Let $h(X)=\sum_{i=0}^n (-1)^{n-i} b_{n-i}X^i$ be the minimal polynomial of $\beta$. Then $\mathrm{Tr}(\beta)=b_1$ and $\mathrm{Tr}(\beta^2)=b_1^2-2b_2$. Hence the assertion follows.
\end{proof}

\begin{lemma}\label{lem-binomials-sum-mod2Equality-general}
Suppose  $1\leq i\leq 2^{s-1}$ and $d_i\in \mathcal{O}_K$. Then for $1\leq m\leq 2^{s-1}(2^f-1)$,
\begin{align}\label{eq-binomials-sum-mod2Equality-general}
&\sum_{k=0}^{v_2(m)}\mathrm{Tr}^K_{\mathbb{Q}_2}\left(u^{\frac{m}{2^k}}\Big(\binom{\frac{m}{2^k}-1}{i-1}d_i+\binom{\frac{m}{2^k}-1}{2i-1}d_i^2\Big)\right)\nn\\
=&\sum_{k=0}^{v_2\big(m+2^{s-1}(2^f-1)\big)}\mathrm{Tr}^K_{\mathbb{Q}_2}\left(u^{\frac{m+2^{s-1}(2^f-1)}{2^k}}\Big(\binom{\frac{m+2^{s-1}(2^f-1)}{2^k}-1}{i-1}d_i+\binom{\frac{m+2^{s-1}(2^f-1)}{2^k}-1}{2i-1}d_i^2\Big)\right)\mod 2.
\end{align}
\end{lemma}
\begin{proof}
By Lemma \ref{lem-integrality-factorial}, we have, for $l\in \mathbb{N}$,
\begin{align*}
-\binom{l-1}{i-1}+\binom{2l-1}{2i-1}=2i\left(-\frac{1}{2}\frac{(l-1)!}{i!(l-i)!}+\frac{(2l-1)!}{(2i)!(2l-2i)!}\right)\equiv 0\mod 2.
\end{align*}
Using this and Lemma \ref{lem-a+a^2-trace-even}, we get
\begin{align*}
&\mathrm{Tr}^K_{\mathbb{Q}_2}\left(u^{\frac{m}{2^k}}\binom{\frac{m}{2^k}-1}{2i-1}d_i^2+u^{\frac{m}{2^{k+1}}}\binom{\frac{m}{2^{k+1}}-1}{i-1}d_i\right)\\
=& \binom{\frac{m}{2^k}-1}{2i-1}\mathrm{Tr}^K_{\mathbb{Q}_2}(u^{\frac{m}{2^k}}d_i^2)+\binom{\frac{m}{2^{k+1}}-1}{i-1}\mathrm{Tr}^K_{\mathbb{Q}_2}(u^{\frac{m}{2^{k+1}}}d_i)\\
=& \binom{\frac{m}{2^k}-1}{2i-1}\mathrm{Tr}^K_{\mathbb{Q}_2}(u^{\frac{m}{2^{k+1}}}d_i)+\binom{\frac{m}{2^{k+1}}-1}{i-1}\mathrm{Tr}^K_{\mathbb{Q}_2}(u^{\frac{m}{2^{k+1}}}d_i)
\equiv 0\mod 2.
\end{align*}
Applying this to $k=0,\dots,v_2(m)-1$, and noting that 
\begin{align*}
\binom{\frac{m}{2^{v_2(m)}}-1}{2i-1}=\frac{2i\cdot 2^{v_2(m)}}{m}\binom{\frac{m}{2^{v_2(m)}}}{2i}
\end{align*}
 is even, we get
\begin{align}\label{eq-binomials-sum-mod2Equality-general-LHS}
\mbox{LHS of (\ref{eq-binomials-sum-mod2Equality-general})}\equiv \binom{m-1}{i-1}\mathrm{Tr}(u^md_i)\mod 2.
\end{align}
The same arguments yield
\begin{align}\label{eq-binomials-sum-mod2Equality-general-RHS}
\mbox{RHS of (\ref{eq-binomials-sum-mod2Equality-general})}\equiv \binom{m+2^{s-1}(2^f-1)-1}{i-1}\mathrm{Tr}(u^{m+2^{s-1}(2^f-1)}d_i)\mod 2.
\end{align}
Using Lemma \ref{lem-binomials-periodicity-general} repeatedly, we get
\begin{align*}
\binom{m+2^{s-1}(2^f-1)-1}{i-1}+\binom{m-1}{i-1}\equiv 0\mod 2.
\end{align*}
Since $u^{m+2^{s-1}(2^f-1)}=u^m$, we obtain $(\ref{eq-binomials-sum-mod2Equality-general-LHS})\equiv(\ref{eq-binomials-sum-mod2Equality-general-RHS})\mod 2$, and hence (\ref{eq-binomials-sum-mod2Equality-general}).
\end{proof}
So we have completed the proof of Theorem \ref{thm-G(m)-essDescendiblePeriodic} in the case $a_n$ is odd.

\subsection{Even \texorpdfstring{$a_n$}{an}}
That $a_n$ is even is equivalent to $\theta\in \mathfrak{P}$.  Applying Theorem \ref{thm-periodicity-G(m)-oddTheta-general} with $b=0$ and $\xi=\theta$, we have
\begin{align*}
\mathrm{Tr}^L_K(\theta^{2m})\in 2^{v_2(m)+2}\mathcal{O}_K,\
\mathrm{Tr}^L_K(\theta^m)\in 2^{v_2(m)+1}\mathcal{O}_K. 
\end{align*}

\begin{theorem}\label{thm-generalExtension-eventualPeriodicity-G(m)-evenTheta}
Suppose $e\leq 2^{s-1}$. Then $G^{L}_K(m)\in 2 \mathcal{O}_K$ is even for $m\geq 2^s+1$.
In particular, $G^L_K(m+2^s)\equiv G^L_K(m)\mod 2 \mathcal{O}_K$ for $m\geq 2^s+1$.
\end{theorem}
\begin{proof}
Expand $\frac{\mathrm{Tr}^L_K(\theta^m)}{2^{v_2(m)+1}}$ as a polynomial of $c_i$'s, $1\leq i\leq e$. By  Lemma \ref{lem-crossingTerms-even-general}, the coefficients of the terms with at least two distinct factors $c_i$ and $c_j$ lie in $2 \mathcal{O}_K$.
By Newton's identity (\ref{eq-Newton-identity-1}), the coefficient of a power of a single $c_i$ in $\frac{\mathrm{Tr}^L_K(\theta^m)}{2^{v_2(m)+1}}$
\begin{comment}
\begin{equation}\label{eq-minus-G(m)-mod2-evenTheta}
  \frac{\mathrm{Tr}(\theta^m)}{2^{v_2(m)+1}}
\end{equation}
\end{comment}
 is 0 if $i\nmid m$, and is $\frac{(-1)^{m}i(-2c_i)^{\frac{m}{i}}}{2^{v_2(m)+1}}$ if $i\mid m$. Note that 
$\frac{m}{i}-v_2(\frac{m}{i})< 2$ exactly when $\frac{m}{i}=1$ or 2. When $m\geq 2^s+1$ there is no such $i$, since $i\leq e\leq 2^{s-1}$ by the assumption. Thus all coefficients of powers of a single $c_i$ in $\frac{\mathrm{Tr}^L_K(\theta^m)}{2^{v_2(m)+1}}$ lie in $2 \mathcal{O}_K$, and the proof is completed.
\end{proof}
\begin{theorem}\label{thm-general-essential-descendiblePeriodicity-G(m)}
Suppose $e\leq 2^{s-1}$. Then for any odd $m\leq 2^s$,
\begin{align*}
\sum_{k=0}^{2^km\leq 2^s}G(2^km)\in 2 \mathcal{O}_K.
\end{align*}
\end{theorem}
\begin{proof}
By Theorem \ref{thm-generalExtension-eventualPeriodicity-G(m)-evenTheta}, the assertion is equivalent to a more convenient form: 
\begin{align*}
\sum_{k=0}^{\infty}G(2^km)\equiv 0\mod 2 \mathbb{Z}_2.
\end{align*}
Recall $G=\mathrm{Tr}^{K}_{\mathbb{Q}_2}\circ G^L_K$, and we expand
\begin{align*}
\sum_{k=0}^{\infty}G^L_K(2^km)
\end{align*}
 as a power series of indeterminants $c_1,\dots,c_{e}$. By  Lemma \ref{lem-crossingTerms-even-general}, we need only consider the powers of every single $c_i$. Then by the proof of Theorem \ref{thm-generalExtension-eventualPeriodicity-G(m)-evenTheta}, we need only to consider such powers with $i=2^km$ for some $k\geq 0$, and for such an $i$, we need only need to consider the power of $c_i$ appearing in $G^L_K(2^km)$ and $G^L_K(2^{k+1}m)$. So we need to consider
\begin{align*}
\frac{(-1)^{2^k m}2^km(-2c_{2^km})}{2^{k+1}}+
\frac{(-1)^{2^{k+1}m}2^km(-2c_{2^km})^2}{2^{k+2}}
\equiv m(c_{2^km}+c_{2^km}^2)\mod 2 \mathcal{O}_K.
\end{align*}
By Lemma \ref{lem-binomials-sum-mod2Equality-general}, $\mathrm{Tr}^K_{\mathbb{Q}_2}(c_{2^km}+c_{2^km}^2)\in 2 \mathbb{Z}_2$, and the proof is completed.
\end{proof}
So we have completed the proof of Theorem \ref{thm-G(m)-essDescendiblePeriodic} in the case $a_n$ is even.

\begin{appendix}
\section{Tables}
In the following tables, we use the notations such as
\begin{align*}
\big\{6k+\{1,5\}\big\}_{k\geq 0}:=\{6k+1: k\geq 0\}\cup \{6k+5: k\geq 0\}.
\end{align*}
For a set $M$, we will describe  a function $h:M\rightarrow \mathbb{Z}/2 \mathbb{Z}$ by displaying
\begin{align*}
\mathrm{Supp}(h):=\{m\in M: h(m)=1 \}.
\end{align*}
By Proposition \ref{prop-F(m,a)=sumOfG(m,a)}, $F(m,\mathbf{a})$ is even if $m$ is odd, so we display only the results for even $m$.

\begin{table}[htbp!]
\centering
\begin{tabular}{|c|c|c|}
\hline
$\mathbf{a}$ & Supp$\{G(m,\mathbf{a})\mod 2: m\in \mathbb{N}\}$ & Supp$\{F(2m,\mathbf{a})\mod 2: m\in \mathbb{N}\}$\\
\hline
$0$ & $\emptyset$  & $\emptyset$ \\
$1$ & $\emptyset$  & $\emptyset$ \\
$2$ & $\{1,2\}$ & $\{2\}$ \\
$3$ & $\{2k+1\}_{k\geq 0}$ & $2 \mathbb{N}$ \\
\hline
\end{tabular}
\caption{$n=1$}
\label{tab:fMod2-n=1}
\end{table}

\begin{table}[htbp!]
\centering
\begin{tabular}{|c|c|c|}
\hline
$\mathbf{a}$ & Supp$\{G(m,\mathbf{a})\mod 2: m\in \mathbb{N}\}$ & Supp$\{F(2m,\mathbf{a})\mod 2: m\in \mathbb{N}\}$\\
\hline
$(0,0)$ & $\emptyset$  & $\emptyset$ \\
$(0,1)$ & $\big\{4k+\{1,2,3\}\big\}_{k\geq 0}$ & $\{4k+2\}_{k\geq 0}$\\
$(0,2)$ & $\{2,4\}$ & $\{4\}$ \\
$(0,3)$ & $\{2k+1\}_{k\geq 0}$ & $2 \mathbb{N}$ \\
$(1,0)$ & $\emptyset$  & $\emptyset$ \\
$(1,1)$ & $\big\{6k+\{1,5\}\big\}_{k\geq 0}$ & $\big\{6k+\{2,4\}\big\}_{k\geq 0}$\\
$(1,2)$ & $\{2\}\cup \{2k+1\}_{k\geq 1}$ & $\{2k\}_{k\geq 2}$ \\
$(1,3)$ & $\big\{6k+\{1,2,3,4\}\big\}_{k\geq 0}$ & $\big\{6k+\{2,6\}\big\}_{k\geq 0}$\\
$(2,0)$ & $\{1,2\}$ & $\{2\}$ \\
$(2,1)$ & $\emptyset$  & $\emptyset$ \\
$(2,2)$ & $\{1,4\}$ & $\{2,4\}$ \\
$(2,3)$ & $\{4k+2\}_{k\geq 0}$ & $4\mathbb{N}$ \\
$(3,0)$ & $\{2k+1\}_{k\geq 0}$ & $2 \mathbb{N}$ \\
$(3,1)$ & $\emptyset$  & $\emptyset$ \\
$(3,2)$ & $\{1,2\}$ & $\{2\}$ \\
$(3,3)$ & $\big\{6k+\{2,3,4,5\}\big\}_{k\geq 0}$ & $\big\{6k+\{4,6\}\big\}_{k\geq 0}$\\
\hline
\end{tabular}
\caption{$n=2$}
\label{tab:fMod2-n=2}
\end{table}

\begin{table}[htbp!]
\centering
\begin{tabular}{|c|c|c|}
\hline
$\mathbf{a}$ & Supp$\{G(m,\mathbf{a})\mod 2: m\in \mathbb{N}\}$ & Supp$\{F(2m,\mathbf{a})\mod 2: m\in \mathbb{N}\}$\\
\hline
$(0,0,0)$ & $\emptyset$  & $\emptyset$ \\
$(0,0,1)$ & $\emptyset$  & $\emptyset$ \\
$(0,0,2)$ & $\{3,6\}$ & $\{6\}$ \\
$(0,0,3)$ & $\big\{6k+3\big\}_{k\geq 0}$ & $\big\{6k\big\}_{k\geq 1}$\\
$(0,1,0)$ & $\big\{4k+\{1,2,3\}\big\}_{k\geq 0}$ & $\{4k+2\}_{k\geq 0}$\\
$(0,1,1)$ & $\big\{14k+\{1,2,3,4,5,6,11,12\}\big\}_{k\geq 0}$ &$\big\{14k+\{2,6,8,10\}\big\}_{k\geq 0}$\\
$(0,1,2)$ & $\{1,2\}$ & $\{2\}$ \\
$(0,1,3)$ & $\big\{14k+\{1,2,4,6,7,11,12,13\}\big\}_{k\geq 0}$ &$\big\{14k+\{2,8,12,14\}\big\}_{k\geq 0}$\\
$(0,2,0)$ & $\{2,4\}$ & $\{4\}$ \\
$(0,2,1)$ & $\big\{6k+\{2,4,5\}\big\}_{k\geq 0}$ & $\big\{6k+4\big\}_{k\geq 0}$\\
$(0,2,2)$ & $\{2,3,4,6\}$ & $\{4,6\}$ \\
$(0,2,3)$ & $\big\{6k+\{2,3,4,5\}\big\}_{k\geq 0}$ & $\big\{6k+\{4,6\}\big\}_{k\geq 0}$\\
$(0,3,0)$ & $\{2k+1\}_{k\geq 0}$ & $2 \mathbb{N}$ \\
$(0,3,1)$ & $\big\{14k+\{1,3,4,8,9,10,12,13\}\big\}_{k\geq 0}$ &$\big\{14k+\{2,4,6,12\}\big\}_{k\geq 0}$\\
$(0,3,2)$ & $\{1\}\cup \{4k+2\}_{k\geq 1}$ & $\{2\}\cup\{4k\}_{k\geq 1}$ \\
$(0,3,3)$ & $\big\{14k+\{1,4,5,7,8,9,10,12\}\big\}_{k\geq 0}$ &$\big\{14k+\{2,4,10,14\}\big\}_{k\geq 0}$\\
$(1,0,0)$ & $\emptyset$  & $\emptyset$ \\
$(1,0,1)$ & $\big\{14k+\{2,3,8,9,10,11,12,13\}\big\}_{k\geq 0}$ &$\big\{14k+\{4,6,8,12\}\big\}_{k\geq 0}$\\
$(1,0,2)$ & $\{4\}\cup \{2k+1\}_{k\geq 1}$ & $\{2k\}_{k\geq 3}$ \\
$(1,0,3)$ & $\big\{14k+\{2,4,5,6,7,9,10,13\}\big\}_{k\geq 0}$ &$\big\{14k+\{4,10,12,14\}\big\}_{k\geq 0}$\\
$(1,1,0)$ & $\big\{6k+\{1,5\}\big\}_{k\geq 0}$ & $\big\{6k+\{2,4\}\big\}_{k\geq 0}$\\
$(1,1,1)$ & $\big\{4k+\{1,2,3\}\big\}_{k\geq 0}$ & $\{4k+2\}_{k\geq 0}$\\
$(1,1,2)$ & $\{1,3,4,5\}\cup \big\{6k+\{2,3,4,5\}\big\}_{k\geq 1}$ & $\{2\}\cup \big\{6k+\{4,6\}\big\}_{k\geq 0}$ \\
$(1,1,3)$ & $\big\{8k+\{1,2,4,5\}\big\}_{k\geq 0}$ & $\big\{8k+\{2,8\}\big\}_{k\geq 0}$\\
$(1,2,0)$ & $\{2\}\cup \{2k+1\}_{k\geq 1}$ & $\{2k\}_{k\geq 2}$ \\
$(1,2,1)$ & $\emptyset$  & $\emptyset$ \\
$(1,2,2)$ & $\{2,4\}$ & $\{4\}$ \\
$(1,2,3)$ & $\big\{14k+\{3,4,5,6,7,8,11,12\}\big\}_{k\geq 0}$ &$\big\{14k+\{6,8,10,14\}\big\}_{k\geq 0}$\\
$(1,3,0)$ & $\big\{6k+\{1,2,3,4\}\big\}_{k\geq 0}$ & $\big\{6k+\{2,6\}\big\}_{k\geq 0}$\\
$(1,3,1)$ & $\big\{8k+\{1,4,5,6\}\big\}_{k\geq 0}$ & $\big\{8k+\{2,4\}\big\}_{k\geq 0}$\\
$(1,3,2)$ & $\{1,2\}$ & $\{2\}$ \\
$(1,3,3)$ & $\{2k+1\}_{k\geq 0}$ & $2 \mathbb{N}$ \\
\hline
\end{tabular}
\caption{$n=3$}
\label{tab:fMod2-n=3}
\end{table}

\begin{table}[htbp!]
\centering
\begin{tabular}{|c|c|c|}
\hline
$\mathbf{a}$ & Supp$\{G(m,\mathbf{a})\mod 2: m\in \mathbb{N}\}$ & Supp$\{F(2m,\mathbf{a})\mod 2: m\in \mathbb{N}\}$\\
\hline
$(2,0,0)$ & $\{1,2\}$ & $\{2\}$ \\
$(2,0,1)$ & $\big\{6k+\{1,2,4\}\big\}_{k\geq 0}$ & $\big\{6k+\{2\}\big\}_{k\geq 0}$\\
$(2,0,2)$ & $\{1,2,3,6\}$ & $\{2,6\}$ \\
$(2,0,3)$ & $\big\{6k+\{1,2,3,4\}\big\}_{k\geq 0}$ & $\big\{6k+\{2,6\}\big\}_{k\geq 0}$\\
$(2,1,0)$ & $\emptyset$  & $\emptyset$ \\
$(2,1,1)$ & $\emptyset$  & $\emptyset$ \\
$(2,1,2)$ & $\{3\}\cup \big\{4k+\{1,2,3\}\big\}_{k\geq 1}$ & $\{4k+2\}_{k\geq 1}$ \\
$(2,1,3)$ & $\big\{14k+\{3,5,7,13\}\big\}_{k\geq 0}$ &$\big\{14k+\{6,10,12,14\}\big\}_{k\geq 0}$\\
$(2,2,0)$ & $\{1,4\}$ & $\{2,4\}$ \\
$(2,2,1)$ & $\big\{6k+\{1,5\}\big\}_{k\geq 0}$ & $\big\{6k+\{2,4\}\big\}_{k\geq 0}$\\
$(2,2,2)$ & $\{1,3,4,6\}$ & $\{2,4,6\}$ \\
$(2,2,3)$ & $\{2k+1\}_{k\geq 0}$ & $2 \mathbb{N}$ \\
$(2,3,0)$ & $\{4k+2\}_{k\geq 0}$ & $4\mathbb{N}$ \\
$(2,3,1)$ & $\big\{14k+\{2,5,6,8,9,10,11,13\}\big\}_{k\geq 0}$ &$\big\{14k+\{4,8,10,12\}\big\}_{k\geq 0}$\\
$(2,3,2)$ & $\{2\}\cup \{2k+1\}_{k\geq 1}$ & $\{2k\}_{k\geq 2}$ \\
$(2,3,3)$ & $\big\{14k+\{2,3,6,7,8,9,10,11\}\big\}_{k\geq 0}$ &$\big\{14k+\{4,6,8,14\}\big\}_{k\geq 0}$\\
$(3,0,0)$ & $\{2k+1\}_{k\geq 0}$ & $2 \mathbb{N}$ \\
$(3,0,1)$ & $\big\{14k+\{1,2,4,5,6,10,11,13\}\big\}_{k\geq 0}$ &$\big\{14k+\{2,8,10,12\}\big\}_{k\geq 0}$\\
$(3,0,2)$ & $\{1,4\}$ & $\{2,4\}$ \\
$(3,0,3)$ & $\big\{14k+\{1,2,3,7,8,10,12,13\}\big\}_{k\geq 0}$ &$\big\{14k+\{2,6,12,14\}\big\}_{k\geq 0}$\\
$(3,1,0)$ & $\emptyset$  & $\emptyset$ \\
$(3,1,1)$ & $\big\{8k+\{2,3,4,7\}\big\}_{k\geq 0}$ & $\big\{8k+\{4,6\}\big\}_{k\geq 0}$\\
$(3,1,2)$ & $\{3,4\}\cup \big\{6k+\{1,2,3,4\}\big\}_{k\geq 1}$ & $\{6\}\cup \big\{6k+\{2,6\}\big\}_{k\geq 1}$ \\
$(3,1,3)$ & $\{4k+2\}_{k\geq 0}$ & $4\mathbb{N}$ \\
$(3,2,0)$ & $\{1,2\}$ & $\{2\}$ \\
$(3,2,1)$ & $\big\{14k+\{1,3,4,5,6,8,9,12\}\big\}_{k\geq 0}$ &$\big\{14k+\{2,4,6,10\}\big\}_{k\geq 0}$\\
$(3,2,2)$ & $\{2,4\}\cup \{2k+1\}_{k\geq 0}$ & $\mathbb{N}\backslash\{4\}$ \\
$(3,2,3)$ & $\big\{14k+\{1,7,9,11\}\big\}_{k\geq 0}$ &$\big\{14k+\{2,4,8,14\}\big\}_{k\geq 0}$\\
$(3,3,0)$ & $\big\{6k+\{2,3,4,5\}\big\}_{k\geq 0}$ & $\big\{6k+\{4,6\}\big\}_{k\geq 0}$\\
$(3,3,1)$ & $\emptyset$  & $\emptyset$ \\
$(3,3,2)$ & $\{2,5\}\cup \big\{6k+\{1,5\}\big\}_{k\geq 1}$ & $\{4\}\cup \big\{6k+\{2,4\}\big\}_{k\geq 1}$ \\
$(3,3,3)$ & $\big\{8k+\{3,4,6,7\}\big\}_{k\geq 0}$ & $\big\{8k+\{6,8\}\big\}_{k\geq 0}$\\
\hline
\end{tabular}
\caption{$n=3$ continued}
\label{tab:fMod2-n=3-continued}
\end{table}
\clearpage

\section{Non-cellularity of elliptic curves}\label{sec:nonSimpleCellular-ellipticCurve}
In this appendix, we will show that  elliptic curves over perfect fields are not cellular in the sense of \cite{MoS20}, as we promised in Example \ref{example-ellipticCurve-dlogRational}. We begin by recalling some necessary definitions.

\begin{definition}[{\cite[Def. 2.9]{MoS20}}]
Let $\mathbf{Sm}_k$ be the category of smooth schemes of finite type over $k$. A Nisnevich sheaf $M$  of abelian groups on $\mathbf{Sm}_k$ is called \emph{strictly $\mathbb{A}^1$-invariant} if the projection induces isomorphisms $H^i_{\mathrm{Nis}}(X,M)\cong H^i_{\mathrm{Nis}}(X\times_k \mathbb{A}^1_k,M)$ for $\forall i\geq 0$. A scheme $X\in \mathbf{Sm}_k$ is \emph{cohomologically trivial} if $H^i_{\mathrm{Nis}}(X,M)=0$ for any strictly $\mathbb{A}^1$-invariant Nisnevich sheaf $M$ and $\forall i>0$.
\end{definition}

\begin{definition}[{\cite[Def. 2.11]{MoS20}}]
\label{def-cellularStructure}
Let $k$ be a field, and $X$ be a smooth $k$-scheme. A \emph{cellular structure} on $X$ consists of an increasing filtration
$\emptyset=\Omega_{-1}\subsetneq \Omega_0\subsetneq\dots\subsetneq \Omega_s=X$
by open subschemes such that for each $i\in \{0,\dots,s\}$, the reduced induced closed subscheme $X_i:=\Omega_i-\Omega_{i-1}$ is $k$-smooth, affine, everywhere of codimension $i$ and cohomologically trivial. A scheme admitting a cellular structure is called a \emph{cellular scheme}.
\end{definition}

\begin{proposition}
An elliptic curve $E$ over a perfect field $k$ is not cellular.
\end{proposition}

We will use \cite[Lemm 2.13]{MoS20}, which states that if a smooth affine scheme $U$ over perfect $k$ is cohomologically trivial, then every vector bundle over $U$ is trivial. 

First we consider the case $k$ is infinite. Every nonempty affine open subset of $E$ is not the Spec of a PID, and thus has a nontrivial line bundle. So we are left to consider the case $k$ is a finite field. In this case,  the affine open subset $E\backslash E(k)$ is the Spec of a PID, and the argument for infinite fields fails. Our idea is to apply \cite[Lemm 2.13]{MoS20} to the self-product of $E\backslash E(k)$.

In the following, let $E$ be an elliptic curve over $\mathbb{F}_q$, $Z$ be a nonempty 0-dimensional subscheme of $E$, and $U=E\backslash Z$.
By \cite[Remark 2.10]{MoS20}, the product of two cohomologically trivial smooth schemes is again cohomologically trivial\footnote{This fact is asserted without proof in \cite[Remark 2.10]{MoS20}. I thank Fabien Morel for explaining the proof to me.}.

\begin{lemma}\label{lem-numChow-EtimesE}
Assume that $U$ is cohomologically trivial. Then
\begin{align*}
\mathrm{CH}_1(E\times E)_{\mathbb{Q}}/\mathrm{Num}= \mathbb{Q}^2.
\end{align*}
\end{lemma}
\begin{proof}
There is a right exact sequence
\begin{align*}
\mathrm{CH}_i(Z\times E)\oplus \mathrm{CH}_i(E\times Z)\rightarrow \mathrm{CH}_i(E\times E)\rightarrow \mathrm{CH}_i(U\times U)\rightarrow 0.
\end{align*}
By the assumption and the above-mentioned result, $U\times U$ is cohomologically trivial. By \cite[Lemma 2.13]{MoS20}, every vector bundle on $U\times U$ is trivial, thus $K_0(U\times U)=\mathbb{Z}$. Then by the Chern character isomorphism, $\mathrm{CH}(U\times U)_{\mathbb{Q}}=\mathbb{Q}=\mathrm{CH}_2(U\times U)_{\mathbb{Q}}$. From the above sequence, $\mathrm{CH}_1(E\times E)_{\mathbb{Q}}$ is generated by the divisors $\{x\}\times E$ and $E\times \{x\}$ for points $x\in Z$. Hence the conclusion follows.
\end{proof}

However, taking the divisors $E\times \{0\}$, and $\{0\}\times E$, and $\Delta_E$ in  $E\times E$, the intersection matrix is nonsingular. It follows that the rank of $\mathrm{CH}_1(E\times E)_{\mathbb{Q}}/\mathrm{Num}\geq 3$ (in fact one can show that this rank is at least 4). 
This leads to a contradiction by Lemma \ref{lem-numChow-EtimesE} if $U$ is cohomologically trivial. Hence we conclude that $E$ is not cellular.

\end{appendix}

\end{document}